\documentclass[11pt]{article}

\usepackage{amssymb,amsmath,latexsym, color}
\usepackage{epsfig}

\newtheorem{theorem}{Theorem}[section]

\newtheorem{lemma}[theorem]{Lemma}
\newtheorem{proposition}[theorem]{Proposition}
\newtheorem{definition}[theorem]{Definition}

\newtheorem{remark}[theorem]{Remark}
\newtheorem{example}[theorem]{Example}
\numberwithin{equation}{section}

\def\proof{{\medskip\noindent {\bf Proof. }}}

\def\qed{{\hfill $\square$ \bigskip}}

\def\square{{\vcenter{\vbox{\hrule height.3pt
        \hbox{\vrule width.3pt height5pt \kern5pt
           \vrule width.3pt}
        \hrule height.3pt}}}}

\def\sA {{\cal A}} \def\sB {{\cal B}} \def\sC {{\cal C}}
\def\sD {{\cal D}} \def\sE {{\cal E}} \def\sF {{\cal F}}
 \def\sH {{\cal H}} 
  \def\sL {{\cal L}}

\def\bD {{\mathbb D}}

 \def\bN {{\mathbb N}} 
\def\bP {{\mathbb P}}  \def\bR {{\mathbb R}}

 \def\bZ {{\mathbb Z}}

\def\wh{\widehat}

\def\E{{\mathbb E}}
\def\P{{\mathbb P}}

\def\bea{\begin{align*}}
\def\eea{\end{align*}}
\def\bee{\begin{equation}}
\def\eee{\end{equation}}

\def\Lip{{\rm Lip}}
\def\1{{\bf 1}}

\def\proof{{\medskip\noindent {\bf Proof. }}}

\def\R{{\mathbb R}}

\def\eps{\varepsilon}

\def\r{\rho}

\def\wh{\widehat}

\def\nn{{\nonumber}}

\makeatletter
\@addtoreset{equation}{section}

\makeatother

\begin{document}
\bibliographystyle{plain}

\title{\bf Discrete Approximation of Symmetric Jump Processes on
Metric Measure Spaces
}
\author{{\bf Zhen-Qing Chen},\thanks{Research partially supported
by NSF Grant  DMS-0906743.} \quad  {\bf Panki Kim}
 \thanks{Research
 supported by Basic Science Research Program through the National Research Foundation of Korea(NRF)
Êgrant funded by the Korea government(MEST)(2010-0001984)} \quad and  \quad {\bf Takashi~Kumagai}
\thanks{Research partially supported by the
Grant-in-Aid for Scientific Research (B) 22340017.}}
\date{(August 31, 2010)}

\maketitle

\begin{abstract}
In this paper we give general criteria on tightness and weak convergence of
discrete Markov chains to symmetric jump processes on metric measure spaces
under mild conditions. As an application, we
investigate discrete approximation for a large class of
symmetric jump processes. We also discuss some application of
our results to the scaling limit of random walk in random conductance.
\end{abstract}
\vspace{.6truein}



\section{Introduction}
For a Hunt process $X$ on $\bR^d$, consider the following question: 
\[\mbox{{\bf (Q1)}~Can $X$ be approximated by a sequence of
of Markov chains $X^{(k)}$ on $k^{-1}\bZ^d$?~~~~~~~~~~~~~~}\]
A closely related question is the following. Let $X^{(k)}$ be a sequence of Markov chains on $k^{-1}\bZ^d$.
\[\mbox{{\bf (Q2)}~When does $X^{(k)}$ converge weakly to a \lq nice\rq\, Hunt process $X$ on $\bR^d$ as $k\to\infty$?  ~~~~~~~}\]
In this paper, we 
address these two questions  when $X$ is a symmetric   process of pure jump.

Let us briefly mention some work 
 on these problems when $X$ is a diffusion.
When $X$ is a diffusion corresponding to an operator in non-divergence form,
these problems were studied, 
for example,
in the book of Stroock-Varadhan (\cite[Chapter 11]{SV}) by solving the corresponding martingale problem.
When $X$ is a symmetric diffusion
corresponding to a
uniformly elliptic divergence form operator, {\bf (Q1)} is
solved completely by Stroock-Zheng \cite{SZ}.
Let $X^{(k)}_t$ be a continuous time symmetric Markov chain on $k^{-1}\bZ^d$
with conductances $\sC^{(k)}(x,y)$;
This means that $X^{(k)}$ stays at a state $x$ for an exponential length of time with parameter
$\sC^{(k)}(x):=\sum_{z\ne x} \sC^{(k)}(x,z)$
and then jumps to the next state $y$ with probability
 $\sC^{(k)}(x,y)/  \sC^{(k)}(x )$.
In \cite{SZ}, they also
answered  {\bf (Q2)} when $\sC^{(k)}(\cdot,\cdot)$ is
 of finite range (i.e. $\sC^{(k)}(x,y)=0$ if $|x-y|\ge R_0/k$ for some $R_0>0$)
 and  has certain uniform regularity.
The core of their paper is to establish a
 discrete version of the De Giorgi-Moser-Nash theory.
Recently, in \cite{BK},
 the main results in \cite{SZ} are
extended
in two ways:
chains with unbounded range were allowed and the strong uniform regularity conditions in \cite{SZ} were weakened.
This was further extended in \cite{BKU} so that the limiting process $X$ had a continuous part and a jump part.
For both \cite{BK,BKU}, a crucial step is
 to obtain a priori estimate of the solution of the heat equation,
which can be derived
 thanks to the recent developments of the De Giorgi-Moser-Nash theory for jump processes.
When $X$ is reflected Brownian motion on a domain, {\bf (Q1)} was solved in \cite{BC}.

Now consider the case where $X$ is a 
 symmetric Hunt process of pure jump. 
Let $(\sE,\sF)$ be its associated symmetric Dirichlet form on
$L^2(\R^d; m)$,
where $m$ is a Radon
measure on $\bR^d$ and
\begin{eqnarray}
\sF&:=& \left\{ u\in L^2(\bR^d, m): \,\int_{\bR^d\times \bR^d \setminus \wh d}{(u(x)-u(y))^2}
J(dx,dy)<\infty \right\},\label{eq:niefbis}  \\
\sE (u, v)&:=& \frac1{2} \int_{\bR^d\times \bR^d \setminus \wh d}
{(u(x)-u(y))(v(x)-v(y))} {J(dx,dy)} \qquad \hbox{for } u, v\in \sF.
\nn
\end{eqnarray}
Here
 $\wh d$ is the diagonal set in $\bR^d \times \bR^d$,
$J(\cdot,\cdot)$ is a measure on $\bR^d\times\bR^d$ such that
$J(A,B)=J(B,A)$.
The paper \cite{KH} considered {\bf (Q1)}--{\bf (Q2)} when
$J(dx,dy)=j(x,y)dxdy$,   $j(x,y)\asymp |x-y|^{-d-\alpha}$ for some $0<\alpha<2$ 
and $m(dx)=dx$. 
(Here and in the following, $f\asymp g$ means that there are $c_1,c_2>0$
so that $c_1g( x)\leq f (x)\leq c_2 g(x)$ in the
common domain of definition for $f$ and $g$.) This
is extended in \cite{BKK} to more general
Dirichlet form $(\sE,\sF)$.
Again, for both \cite{BKK, KH}, the crucial point
 is to obtain a priori H\"older estimate of the solution of
the heat equation.
However for general symmetric Markov processes, obtaining good a priori
estimate for their transition densities is impossible.
Indeed, even in the case $c_1|x-y|^{-d-\alpha_1}\le j(x,y)\le c_2|x-y|^{-d-\alpha_2}$
for $|x-y|<1$ where $\alpha_1<\alpha_2$, one can construct an example where there
is a bounded harmonic function that is not continuous (see \cite[Theorem 1.9]{BBCK}). \\

In this paper, we will answer {\bf (Q1)} affirmatively for a very general
class of symmetric Markov processes whose associated Dirichlet forms
are of the form \eqref{eq:niefbis} (see Theorem \ref{t:latappr}),
and give answer
to {\bf (Q2)} when $X^{(k)}$ and $X$ satisfy
conditions {\bf (A1)}--{\bf (A4)} in Section \ref{tightsec}--\ref{sec4}
(see Theorem \ref{t:wc}).
Our approach does not
rely on the a priori estimate of the heat kernel, instead we adapt  the ideas
 of \cite{BC} and use the Lyons-Zheng decomposition
to obtain tightness (Proposition \ref{tight_aa}).
The drawback is we can only obtain tightness when
the initial distribution is
absolutely continuous with respect to the reference measure.
Note that when we have a priori estimate of the heat kernel (such as examples discussed in
\cite{BKK, KH}), we can obtain tightness for any initial distributions.
To show
finite dimensional distribution convergence, we establish
 the Mosco convergence, which is
equivalent to strong convergence of the semigroups
(Theorems \ref{Mosco1} and \ref{Mosco1--0}). We will obtain these results on
a large class of
 metric measure spaces with volume doubling property.

It is quite important and useful if we can obtain {\bf (Q2)}
 in such a way that
 is applicable to prove convergence
of Markov chains on some random media. In order to establish such results, we need to relax the assumption
for $X^{(k)}$. In Theorem \ref{Mosco1--0}, we prove the Mosco convergence under a milder
 condition on $X^{(k)}$ and a stronger condition on $X$.
Then the following example can be handled.
Let $\{\xi_{xy}\}_{x,y\in
\bZ^d,
x\ne y}$ be
i.i.d. on
a probability space $(\Omega, \sA, {\bf P})$ such that
 $0\le\xi_{x,y}$,
${\bf E}[\xi_{x,y}] =1$ and $\mbox{Var }(\xi_{x,y})<\infty$.
Let $d\ge 2$, $0<\alpha<2$ and
\[
\sC(x,y)=\xi_{xy}|x-y|^{-d-\alpha},\qquad  x,y\in \bZ^d\]
be the random conductance.
Let $X^{(1)}$ be the corresponding Markov chain on $\bZ^d$
 with this conductance.
 Then we can prove that
$X^{(k)}_t=k^{-1}X^{(1)}_{k^\alpha t}$ converges weakly to (a constant time change of)
 symmetric $\alpha$-stable process on $\bR^d$
 equipped with convergence-in-measure topology  ${\bf P}$-a.s.
(see  Theorem \ref{arrc} (i)).
Moreover, if we further assume that $0\le\xi_{xy}\le C_1$ for some deterministic constant $C_1>0$,
we can prove that
$X^{(k)}_t$ converges weakly 
on $\bD([0, 1]; \bR^d)$ equipped with the Skorohod topology
to symmetric $\alpha$-stable process on $\bR^d$ ${\bf P}$-a.s.
(see  Theorem \ref{arrc}(ii)).

The rest of the paper is organized as follows.
In Section 2,
 we present our framework of the base metric measure space $M$ and discuss its graph approximation.
In Section 3, we give a family of Markov chains $X^{(k)}$ on the approximated lattices and give tightness criteria.
In Section 4, we give a symmetric pure jump process $X$ on $M$ and give sufficient condition for finite dimensional distribution convergence of $X^{(k)}$ to $X$.
Section \ref{Sect5} is for our main theorems on weak convergence and discrete approximation of $X$.
In Section \ref{Sect6}, we give tightness and weak convergence of $X^{(k)}$ under the convergence-in-measure topology
which is a topology weaker than the Skorohod topology.
In Section 7, we apply our results to random walk in random conductance. Finally in Appendix we give a
full proof of generalized
Mosco convergence.

For technical convenience, we will often consider stochastic
processes whose initial distribution is a finite measure, not
necessarily normalized to have total mass 1,  for example,
$\varphi(x)m(dx)$ where $\varphi$ is bounded function with compact
support. Translating our results to the usual probabilistic setting
is straightforward and so it is left to the reader.

Throughout paper, we use ``$:=$" to denote a definition, which is  read
as ``is defined to be". The letter $c$, with or without subscripts,
signifies a constant whose value is unimportant and which may change
from location to location, even within a line.
For a metric space
$M$, we use   $C(M)$ to denote the  space of continuous functions
on $M$ and $\Lip(M)$ the space of Lipschitz
continuous functions on $M$. For any collection of numerical functions $\sH$, $\sH^+$ denotes the set of nonnegative functions in $\sH$,   $\sH_b$ denotes the set of bounded functions in $\sH$ and $\sH_c$ denotes the set of functions in $\sH$  with compact support. Moreover, we denote $\sH^+_c:=\sH^+ \cap \sH_c$ and $\sH^+_b:=\sH^+ \cap \sH_b$.
We will use
$\#S$ is the cardinality of a set $S$.

\medskip

\section{Discrete approximation of the space}\label{framedisc}

Let $(M,\rho, m)$ be a metric measure space,  where $(M,\r)$ is a
locally compact separable
connected
metric space and $m$ is a Radon
measure on $M$ with $V(x,r):= m(B(x,r))\in (0,\infty)$ and $m(\partial
B(x,r))=0$ for each $r>0$, $x\in M$. Here and in the sequel,
$B(x,r)$ denotes the open ball of radius $r$ centered at $x$, and
$\partial B(x,r)=\overline{B(x,r)}\setminus B(x,r)$. We assume the
following:

\begin{description}
\item{\bf (MMS.1)} The closure of $B(x,r)$ is compact for every $x\in M$ and $r>0$.

\item{\bf (MMS.2)} $\r$ is geodesic, that is, for any two points
$x,y\in M$, there exists a continuous map $\gamma:[0,\r (x,y)]\to M$
such that
$\gamma (0)=x$, $\gamma (\r (x, y))=y$ and
 $\r(\gamma (s),\gamma (t))=t-s$ for all $0\le s\le t\le
\r(x,y)$.

\item{\bf (MMS.3)} $(M,\rho, m)$ satisfies volume doubling
 property (VD for short), that is,
$$
\text{there is a constant } c_1>0 \text{ such that } V(x,2r)\le c_1V(x,r)\quad  \hbox{for every }  x\in M \hbox{ and }  r>0.
$$
\end{description}

Fix some $x_0\in M$. Condition {\bf (MMS.3)} in particular implies that
$$V(x_0, 2^n)\leq c_1^n V(x_0, 1)= (2^n)^{\log_2 c_1} V(x_0, 1) \qquad
\hbox{for every } n\geq 1.
$$
So  there are constants
$c_0=c_0(x_0)>0$ and $d_0>0$ such that
\begin{equation}\label{e:2.1}
V(x_0, r)\leq c_0 r^{d_0} \qquad \hbox{for every } r\geq 1.
\end{equation}
It follows then
\begin{eqnarray}
\int_M e^{-\lambda \rho (x, x_0)} m(dx)&=& \int_0^\infty e^{-\lambda r} d ( V(B(x_0, r))
=\lambda \int_0^\infty  V (B(x_0, r))\, e^{-\lambda r}\, dr  \nonumber \\ \label{e:2.2}
&\leq &   c\, \lambda \left(1+\int_1^\infty r^{d_0} \, e^{-\lambda r}   dr\right) <\infty.
\end{eqnarray}

Consider approximating graphs $\{(V_k,B_k),\,k\in \bN\}$ of $M$ with
the graph distance $\rho_k$ and the associated partition
$\{U_k(x), x\in V_k; \, k\in \bN\}$ that satisfies the following
properties.
Here $V_k$ is the set of vertices and $B_k$ is the set of edges of
the graph $(V_k, B_k)$.

\begin{description}
\item{\bf (AG.1)}  $(V_k,B_k)$ is connected and has uniformly bounded degree.

\item{\bf (AG.2)} $V_k\subset M$, $\cup_k V_k$ is dense in $M$ and
\begin{equation}\label{eq:equidis22}
\frac{C_1}k\rho_k(x,y)\le \rho(x,y)\le \frac{C_2}k\rho_k(x,y)\qquad \hbox{for every }  x,y\in V_k.
\end{equation}

\item{\bf (AG.3)} $\cup_{x\in V_k} U_k(x)=M$, $m(U_k(x)\cap U_k(y))=0$ for $x\ne y$, and
\begin{equation}\label{eq:equidis25}
\sup\{\rho(\xi,\eta): \xi,\eta\in U_k(x)\}\le C_3/k.
\end{equation}
Moreover, for each $x\in V_k$, $V_k\cap  \mbox{Int}\, U_k(x)
=\{x\}$, and we have
\begin{equation}\label{eq:equidis27}
C_4 \, m(U_k(x))\le V(x,1/k)\le C_5\,  m(U_k(x)).
\end{equation}
\end{description}

\medskip

\begin{theorem}\label{T:2.1}
Suppose $(M,\rho, m)$ is a   metric measure space  satisfying
conditions {\bf (MMS.1)}--{\bf (MMS.3)}. Then $M$ admits
approximating graphs $\{(V_k, B_k), k\geq 1\}$ and
associated
partitions $\{U_k(x), x\in V_k; k\geq 1\}$ satisfying {\bf
(AG.1)}--{\bf (AG.3)}.

\end{theorem}

To prove this theorem, we need the following \lq nice\rq\, open
covering of $M$ (see, for example \cite[Lemma 3.1]{KumS}, for a
proof).

\begin{lemma}\label{ksvdc}
Suppose $(M,\rho, m)$ is a   metric measure space  satisfying
conditions {\bf (MMS.1)}--{\bf (MMS.3)}. Then  there
   exist integers  $N_0, L_0\geq 1$ that depend only on the constant $c_1$ in {\bf (MMS.3)}
such that for each $r>0$ there exists an open covering $\{B(x_i,r),
\, i\geq 1\}$  of $M$ with the following property:
\begin{itemize}
\item No point in $M$ is contained in more than $N_0$ of the
balls $\{B(x_i,r), \, i\in \bN\}$.
\item $\{B(x_i,r/2),
 \, i\in \bN\}$ are disjoint.
\item For each $x\in M$, the number of balls $B(x_i,r)$
which intersects with $B(x,2r)$ is bounded by $L_0$.
\end{itemize}
\end{lemma}

\noindent{\bf Proof of Theorem \ref{T:2.1}.} Let $V^{(r)}=\{x_i,
i\geq 1\}$, where $\{x_i, i\geq 1\}$ are given in Lemma \ref{ksvdc}.
We say two distinct $x, y\in V^{(r)}$ are connected by a {\em  bond}
(which we will denote as $\{x,y\}\in B^{(r)}$)  if $\rho (x,y)\le
3r$. In this way, we can define a graph $(V^{(r)},B^{(r)})$ of
bounded degree. We also define $\{U_{(r)}(x)\}_{x\in V^{(r)}}$, an
associated partition of $M$, as follows;
$U_{(r)}(x_1)=\overline{B(x_1,r)}$ and $U_{(r)}(x_k)=
\overline{B(x_k,r)}\setminus \cup_{i=1}^{k-1}B(x_i,r)$ for $k\ge 2$.
Clearly, $c_1V(x_i,r)\le m(U_{(r)}(x_i))\le V(x_i,r)$ and
$U_{(r)}(x_i)\cap U_{(r)}(x_j)\subset \cup_{k=1}^j\partial B(x_k,r)$
for $i<j$. The definition of $(V^{(r)},B^{(r)})$ and partition
$\{U_{(r)}(x), \, x\in V^{(r)} \}$ depends on the choice of the open
covering of $M$ (and its labeling). In the following, for each
$r>0$, we choose one open covering with the above mentioned property
and fix the graph $(V^{(r)},B^{(r)})$ and a partition $\{U_{(r)}(x),
\, x\in V^{(r)}\}$. For each sequence $(r_m)$ which converges to
zero, the set $\cup_m V^{(r_m)}$ is dense in $M$. Note that since
$\rho$ is geodesic, for each $x\in V^{(r)}$, there exists $y\in
V^{(r)}\setminus \{x\}$ such that $y\in B(x,2r)$. So
$(V^{(r)},B^{(r)})$ is connected. Further, $(V^{(r)},B^{(r)})$ has
bounded degree, i.e. $\sup_{x\in V^{(r)}}\sharp \{y\in V^{(r)}:
\{x,y\}\in B^{(r)}\}<\infty$. Let $\rho^{(r)}$ be the graph distance
of $(V^{(r)},B^{(r)})$; then
\begin{equation}\label{eq:equidis}
\frac{r}2 \, \rho^{(r)}(x,y)\le \rho(x,y)\le 3r\rho^{(r)}(x,y)\qquad \hbox{for }
 x,y\in V^{(r)}.
\end{equation}
 Clearly, this holds if
  $\{x, y\}\in B^{(r)}$.
 In general, the second inequality of (\ref{eq:equidis}) clearly holds
and the first inequality can be verified
as follows.  Let $\gamma$ be a geodesic connecting $x$ and $y$. Set
$k=[1+ r^{-1} \rho (x, y)]$, the largest integer not exceeding $1+
r^{-1} \rho (x, y)$. Let $\{y_i, 0\leq i\leq k\}$ be equally spaced
points on $\gamma$ so that $\rho (y_{i-1}, y_i)=\rho (x, y)/k <r$
for $k=1, \cdots, k$ with $y_0=x$ and $y_k=y$. For each $1\leq i\leq
k=1$, there is some $x_i\in V^{(r)}$ so that $y_i\in B(x_i, r_i)$
(we take $x_0=y_0=x$ and $x_k=y_k=y$). By
the triangle inequality,
$$ \rho (x_{i-1}, x_i)\leq \rho (x_{i-1}, y_{i-1})+\rho (y_{i-1}, y_i)+
\rho (y_i, x_i)   <3r \qquad \hbox{for } i=1, \cdots, k.
$$
This shows that $\rho^{(r)} (x, y) \leq k  \leq 2\rho (x, y)/r$,
establishing the first inequality in \eqref{eq:equidis}.
Let $V_k:=V^{(1/k)}, B_k:=B^{(1/k)}, \rho_k:=\rho^{(1/k)}$ and $U_k(x):=U_{(1/k)}(x)$.
 It is now easy to verify that $(V_k,B_k,\rho_k)$
together with $\{U_k(x), \, x\in V_k\}$
satisfies {\bf (AG.1)}--{\bf (AG.3)}.
\qed

\section{Tightness}\label{tightsec}

For the remainder of this paper, we assume that
 $(M,\rho, m)$ is a   metric measure space  satisfying
conditions {\bf (MMS.1)}--{\bf (MMS.3)} and that
 $\{(V_k, B_k), k\geq 1\}$ are approximating graphs with associated
partitions $\{U_k(x), x\in V_k; k\geq 1\}$ satisfying {\bf
(AG.1)}--{\bf (AG.3)}.

 Let $m_k$ be the measure defined on $V_k$ by
 \begin{equation}\label{e:m_k}
  m_k(A)=\sum_{y\in A}m (U_k(y)) \qquad \hbox{for } A\subset V_k.
 \end{equation}
 For $y\in V_k$, $m_k (\{y\})$ will simply be denoted by $m_k(y)$.

For $k\in \bN$, let $\{ \sC^{(k)}(x,y), \, x,y\in V_k\}$ be a family
of conductance defined on the graph $(V_k, B_k)$; that is,
$\sC^{(k)}(x,y)=\sC^{(k)}(y,x)\ge 0$ for $x, y\in V_k$.
Note that in contrast with  notations in some literatures on graphs
here the set $B_k$ of edges only gives the topological structure of
the graph and has nothing to do with the conductances; that is,
$B_k$ can be different from  $\{(x, y): \, \sC^{(k)}(x,y)>0\}$.
Note also that the graph with vertices $V_k$ and bonds $\{(x, y): \, \sC^{(k)}(x,y)>0\}$
could be disconnected.
 We consider the following quadratic form $(\sE^{(k)},\sF^{(k)})$:
\begin{eqnarray}
\sF^{(k)}&:=& \Big\{ u\in L^2(V_k; m_k)\,;\, \sum_{x,y \in V_k}
(u(x)-u(y))^2 \sC^{(k)}(x,y) m_k(x)m_k(y) < \infty \Big\}
\label{e:Dfk1} \\
\sE^{(k)} (u, v)&:=& \frac1{2}  \sum_{x,y \in V_k}
{(u(x)-u(y))(v(x)-v(y))} \sC^{(k)}(x,y) m_k(x)m_k(y) \quad \hbox{for
} u, v\in \sF^{(k)}. \nn
\end{eqnarray}
It is easy to check by using Fatou's lemma that $(\sF^{(k)},
\sE^{(k)})$ is a Dirichlet form on $L^2(V_k; m_k)$.

\begin{theorem}\label{T:2.4}
Suppose that for each $k\in \bN$ and each compact (or equivalently, finite) set $K\subset V_k$,
\begin{equation}\label{e:2.10}
 \sup_{x\in K} \sum_{y\in V_k} \sC^{(k)}(x,y)
  m_{k}(y) <\infty .
 \end{equation}
Then $C_c(V_k)\subset \sF^{(k)}$ and $(\sE^{(k)},\sF^{(k)})$ is a
regular Dirichlet form on $L^2(V_k; m_k)$. If
\begin{equation}\label{e:cons1}
\sup_{x\in V_k} \sum_{y\in V_k} \sC^{(k)}(x, y)m_k(y)<\infty,
\end{equation}
then the symmetric Hunt process $X^{(k)}$ on $V_k$ associated with
the regular Dirichlet form  $(\sE^{(k)},\sF^{(k)})$ is conservative.
\end{theorem}

\proof    For  $f\in C_c(V_k)$, let $K$ denote
its support. Then under condition \eqref{e:2.10},
 \begin{eqnarray*}
  \sE^{(k)}(f, f)&=& \frac12 \sum_{x, y\in K} (f(x)-f(y))^2 \sC^{(k)}(x, y)
  m_k(x)m_k(y)  \\
  && +  \sum_{x\in K} f(x)^2 \left(\sum_{y\in K^c} \sC^{(k)}(x, y)   m_k(y) \right)m_k(x)\\
    &\leq & 3
  \|f\|_\infty^2 \sum_{x\in K} \left( \sum_{y\in V_k}\sC^{(k)}(x, y)m_k(y)\right)
   m_k(x) <\infty.
  \end{eqnarray*}
  This shows that $f\in \sF^{(k)}$ and so $C_c(V_k)\subset \sF^{(k)}$.
  Let $K_j$ be an increasing sequence of compact (or
equivalently, finite) subsets of $V_k$ with $\cup_{j\geq 1}K_j=V_k$.
For every $u\in \sF^{(k)}_b$,
define $u_j=u-((-1/j)\vee u))\wedge (1/j)$.
By \cite[Theorem 1.4.2(iv)]{FOT}, $u_j$ is $\sE^{(k)}_1$-convergent to $u$.
Since $u\in L^2(V_k; m_k)$,
${\rm supp}[u_j]\subset \{x\in V_k: |u(x)|>1/j\}$   is a finite set.
Consequently $u_j\in C_c(V_k)$ and so
  $(\sE^{(k)}, \sF^{(k)})$ is a regular Dirichlet form
on $L^2(V_k; m_k)$. Thus there is an associated $m_k$-symmetric Hunt
process $X^{(k)}$ on $V_k$.

Fix some $x_0\in V_k$. Note that for $r>0$, by {\bf (AG.1)}--{\bf (AG.3)}
and \eqref{e:2.1}
\begin{eqnarray*}
m_k (B(x_0, r))&:=& \sum_{y\in V_k, \rho_k(x_0, y)\leq r} m_k(y)
\,=\,m\Big( \bigcup_{y\in V_k, \rho_k (x_0, y)\leq r} U_k (y) \Big)\\
&\leq & m \left(B(x_0, C_2r+C_3 ) \right) \,\leq\, c (r+1)^{d_0}.
\end{eqnarray*}
Thus for every $\lambda >0$,
\begin{eqnarray*}
\int_{V_k} e^{-\lambda \rho_k (x, x_0)} m_k(dx)&=& \int_0^\infty e^{-\lambda r} d ( m_k(B(x_0, r))
\,=\,\lambda \int_0^\infty  m_k (B(x_0, r))\, e^{-\lambda r}\, dr  \nonumber \\
&\leq &   c\, \lambda \left( \int_0^\infty (1+r)^{d_0} \, e^{-\lambda r}   dr\right) \,<\,\infty.
\end{eqnarray*}
Note that $\rho_k$ is a discrete metric on $V_k$ and so condition
\eqref{e:cons1} is equivalent to having
$$ \sup_{x\in V_k} \sum_{y\in V_k} (\rho_k (x, y)^2 \wedge 1) \sC^{(k)}(x,
y)m_k(y)<\infty.
$$
 So we conclude from \cite[Theorem
3.1]{MU} that under the condition  \eqref{e:cons1} that $X^{(k)}$ is
conservative. \qed

For notational convenience, fix some $x_0\in M$ and, for $r>0$,
denote $B(x_0, r)$ by $B_r$. Note that by assumption {\bf (MMS.1)},
$\overline B_r$ is compact for every $r>0$.

 Consider  following condition:

\medskip

\noindent{\bf (A1)}. There is $k_0\geq 1$ so that
 for every integer $j\geq 1$,
\begin{equation}\label{e:A1a}
\sup_{k\geq k_0}\sup_{x \in \overline B_j \cap V_k} \sum_{y\in V_k} \sC^{(k)}(x,y) \Big(
\frac {\rho_k(x,y)}k  \wedge 1\Big)^2 m_{k}(y) <\infty
\end{equation}
 and
\begin{equation}\label{e:A1b}
  \sup_{k\geq k_0} \, \sup_{x\in B_{j+2}^c\cap V_k} \sum_{y\in  B_j \cap V_k}
    \sC^{(k)} (x, y)m_k(y)<\infty.
\end{equation}

For every positive function $\varphi \in C_c(M)$,
 we define measures
$$
\P^{(k)}_{\varphi} (\,\cdot\,) \,:= \,\sum_{x\in V_k} \P^{(k)}_x (\,\cdot\,) \varphi(x)m_k(x)
\qquad\text{and}\qquad
\P_{\varphi} (\,\cdot\,) \,:= \,\int_{M} \P_x (\,\cdot\,) \varphi(x)m(dx).
$$

\begin{lemma}\label{l:ct}
Assume condition
{\bf (A1)}
holds. Then for every $g\in \Lip_c(M)$,
there exists a positive constant $c$ such that for every  $k\geq k_0$
and $0 \le s < t<\infty$,
$$
\int_s^t \sum_{y \in V_k} (g(X^{(k)}_u)-g(y))^2 \sC^{(k)}(X^{(k)}_u, y)m_k(y)du
\, \le\, c (t-s).
$$
\end{lemma}

\proof
Let $\Lambda$ be the Lipschitz constant of $g$.
There is an integer $j\geq 1$ so that the topological support $K$ of $g$
is contained in ball $B_j$ centered at $x_0$ with radius $j$.
Let $K_1:=\overline B_{j+1}$ and
$K_2:=\overline B_{j+3}$.
 By
\eqref{eq:equidis22} and \eqref{e:A1a}--\eqref{e:A1b},
\begin{eqnarray*}
 && \sup_{x\in V_k} \sum_{y \in V_k} (g(x)-g(y))^2
\sC^{(k)}(x, y) m_k(y)  \\
&= & \sup_{x\in V_k} \left( \sum_{y \in K_1^c \cap V_k}  g(x)^2 \sC^{(k)}(x, y)m_k(y)+
\sum_{y\in K_1\cap V_k} (g(x)-g(y))^2
\sC^{(k)}(x, y) m_k(y)\right)  \\
&\leq & \|g\|_\infty^2 \sup_{x\in K\cap V_k}  \sum_{y \in K_1^c \cap V_k}    \sC^{(k)}(x, y)m_k(y)
  + \sup_{x\in K_2^c \cap V_k}
\sum_{y\in K_1\cap V_k}  g(y)^2
\sC^{(k)}(x, y) m_k(y)   \\
&&
+\sup_{x\in K_2 \cap V_k} \sum_{y\in K_1\cap V_k}
\left( \Lambda^2\rho (x, y)^2 \wedge 4\|g\|_\infty^2\right) \sC^{(k)}(x, y) m_k(y)\\
  &\leq & c_1  \|g\|_\infty^2
 +c_2 \sup_{x\in K_2\cap  V_k} \sum_{y\in   V_k}
\left( \frac{\rho_k (x, y)  }{k } \wedge 1\right)^2 \sC^{(k)}(x, y) m_k(y)
\leq  c_3,
\end{eqnarray*}
where $c_1$, $c_2$
and $c_3$ are positive constants  independent of $k\geq
k_0$. The conclusion of the lemma follows directly from the above
inequality. \qed

Let $M_\partial:=M \cup \{\partial \} $ be the
 one-point compactification of $M$, and let
\begin{eqnarray*}
 \bD_{M_\partial }[0,\infty)&:=&\big\{f: [0,\infty)\to M_\partial \, \big| \,   f \hbox{ is right continuous having left limits.}   \big\}.
\end{eqnarray*}
 Clearly $X^{(k)}\in \bD_{M_\partial }[0,\infty)$.

Since $\Lip^+_{c} (M)=\{f \in \Lip_{c} (M): f \ge 0\}$ separates
points of $M$, using Stone-Weierstrass theorem, it is easy
to check that $\Lip^+_{c} (M)$ is a dense subset of
$C^+_\infty(M)$ (space of non-negative continuous functions on $M$
on $M$ that vanishes at infinity).

\begin{proposition}\label{tight_aa} Assume {\bf (A1)} holds and let $\zeta^{(k)}$ denote the lifetime
 of the process $X^{(k)}$.
Then, for any 
$\varphi \in C^+_c(M)$, $T>0$, $m \ge 1$ and $\{g_1 ,
\cdots, g_m \} \subset \Lip^+_{c} (M)$, the laws of $\left\{( g_1 , \cdots,
g_m)(X^{(k)})\right\}_{k \ge 1}$ on $\{\zeta^{(k)}>T\}$
with initial distribution $\varphi
(x) m_k (dx)$ is tight in 
$\bD_{\R^m}[0,T]$
equipped with the Skorohod topology. 
Moreover,  the laws of $\big\{X^{(k)}_t, t\in [0, T] \big\}$ on $\{\zeta^{(k)}>T\}$ with initial distribution $\varphi
(x) m_k (dx)$ is tight in $\bD_{M_\partial }[0,T]$ equipped with the 
Skorohod topology. 
\end{proposition}

\proof Without loss of generality, we assume that $m=1$, $T=1$ and $g=g_1$. 
We first show that 
$\left\{(g(X^{(k)}), \P^{(k)}_\varphi); \, k\geq 1\right\} $
 is relatively compact  
 in  $\bD_{\R }[0,T]$ equipped with the Skorohod topology.

Given $t>0$ and a path $\omega \in \bD_{M} [0, 1]$, the time reversal operator $r_t$
is defined by
\begin{equation*}
r_t(\omega)(s):= \left\{\begin{array}{ll} \omega((t-s)-),
& \mbox{if }~0\le s\le t, \\
\omega(0) & \mbox{if}~ s \ge t.
\end{array}\right.\end{equation*}
Here for $r>0$, $\omega(r-):= \lim_{s \uparrow r} \omega(s)$ is the left limit at $r$ and we use the convention that
$\omega(0-):=\omega(0)$

 Since $f|_{V_k} \in \sF^{(k)}$ for every $f \in \Lip_{c} (M)$, by
 the same argument as that for (2.3) in \cite{BC2} (see also \cite{CFKZ}),
  we have  the following forward-backward martingale
decomposition of $f(X^{(k)}_t)$ for every $f \in \Lip_{c} (M)$;
There exists a martingale $M^{k,f}$ such that  on $\{\zeta^{(k)}>T\}$,
\bee\label{e:fbmd}
f(X^{(k)}_t)-f(X^{(k)}_0)=\frac12M^{k,f}_t-\frac12(M^{k,f}_1-M^{k,f}_{(1-t)-})
\circ r_1, \quad t \in [0,1]
 \eee
  By Proposition 2.8 in \cite{CFKZ},
for each $M^{k,f}$,  there exists the continuous predictable
quadratic variation process $\langle  M^{k,f} \rangle_t$. Note that
(for example, see the page 214 of \cite{FOT})
$$
\langle  M^{k,f} \rangle_t-\langle  M^{k,f} \rangle_s=
\int_s^t \sum_{y \in V_k} (f(X^{(k)}_u)-f(y))^2 \sC^{(k)}(X^{(k)}_u, y)m_k(y)du.
$$
Thus by Lemma \ref{l:ct} and \cite[Proposition VI.3.26]{JS},
$\{\langle  M^{k,f} \rangle_t\}_{k \ge 1}$ is $C$-tight in $\bD_{\bR}
[0, 1]$ equipped with the Skorohod topology. As $m_k$ converges
weakly to $m$, by \cite[Theorem VI.4.13]{JS} the laws of $\{ M^{k,f} \}_{k \ge
1}$   is tight in
 $\bD_{\R} [0, 1]$ with the initial distribution $\P^{(k)}_h$
for every $h  \in \Lip^+_{c} (M)$. Thus the laws of $\{ M^{k,f},
\mu^{(k)}_{h_1, h_2} \}_{k \ge 1}$   is tight in the sense of
Skorohod topology on $\bD_{\bR} [0, 1]$ for every $h_1, h_2  \in
\Lip^+_{c} (M)$ where
$$\mu^{(k)}_{h_1, h_2}(A):= \E\left[h_1(X^{(k)}_0(\omega))\1_A(\omega) h_2(X^{(k)}_1(\omega)); \,
  \zeta^{(k)}>1\right],\quad \forall A
\in \sB(\bD_{M} [0, 1]).$$
Note that for every $A \in \sB(\bD_{M} [0, 1])$,
\begin{eqnarray*}
\mu^{(k)}_{h_1, h_2}(A \circ r_1 ) &=&
\E\left[h_1(X^{(k)}_0(\omega))\1_A\circ r_1
(\omega)h_2(X^{(k)}_1(\omega)); \,  \zeta^{(k)}>1\right]\\
&\,=\,& \E\left[h_2(X^{(k)}_0(\omega))\1_A
(\omega)h_1(X^{(k)}_1(\omega)); \, \zeta^{(k)}>1\right]\\
&\, =\,& \mu^{(k)}_{h_2, h_1}(A).
\end{eqnarray*}
 Thus the laws of $\{ M^{k,f},
\mu^{(k)}_{h_1, h_2} \}_{k \ge 1}$  is the same as the laws of $\{
M^{k,f}\circ r_1, \mu^{(k)}_{h_2, h_1} \}_{k \ge 1}$ and so the laws
of $\{ M^{k,f} \circ r_1, \mu^{(k)}_{h_1, h_2} \}_{k \ge 1}$ is
tight in the sense of Skorohod topology on $\bD_{\bR} [0, 1]$ for
every $h_1, h_2 \in \Lip_{c} (M)$, too. So the laws of $\big\{
M^{k,f}, \mu^{(k)}_{\varphi, f} \big\}_{k \ge 1}$ and the laws of
$\big\{ M^{k,f} \circ r_1, \mu^{(k)}_{\varphi, f} \big\}_{k \ge 1}$
are tight. Since the laws of $\big\{g(X^{(k)}), \P^{(k)}_\varphi
\big\}_{k \ge 1}$  restricted to
$\{\zeta^{(k)}>1\}$ are the same as $\big\{g(X^{(k)}),
\mu^{(k)}_{\varphi, g} \big\}_{k \ge 1}$ in $\bD_{\bR} [0, 1]$,  by
\eqref{e:fbmd} $\big\{g(X^{(k)}), \, \P^{(k)}_\varphi\big\}_{k \ge
1}$  restricted to
$\{\zeta^{(k)}>1\}$ is tight (and so relatively compact)
 in the sense of Skorohod topology on 
 $\bD_{\bR} [0, 1]$.
 
 Since  $M_\partial$ is compact and 
the linear span of $\Lip^+_c(M)$ and constants is a dense subset
in $C(M_\partial)$ equipped with uniform topology,
we conclude from \cite[Theorem 3.9.1 and Corollary 3.9.3]{EK}
that the laws of $\big\{X^{(k)}_t, t\in [0, T] \big\}$ on $\{\zeta^{(k)}>T\}$ with initial distribution $\varphi
(x) m_k (dx)$ is tight in $\bD_{M_\partial }[0,T]$ equipped with the
Skorohod topology.
\qed

\section{Semigroup convergence}\label{sec4}

In this section, we discuss semigroup convergence in two settings.

In Section \ref{Sect5},
we will show  that $X^{(k)}$ converges to $X$ in
the sense of finite dimensional distributions. One way to establish
this is to show that corresponding Dirichlet form converges in the
sense of Mosco, a concept introduced in \cite{M}. In \cite{M}, a
symmetric bilinear form $a(u,u)$ defined on a linear subspace
$\sD[a]$ of a Hilbert space $\sH$ is extended to the whole space
$\sH$ by defining $a(u,u) = \infty$ for every $u \in \sH \setminus
\sD[a]$. We will use this extension throughout this paper. In
\cite{M}, Mosco showed that the Mosco convergence of a sequence of
densely defined symmetric closed forms defined on  the same Hilbert
space is equivalent to the convergence of the sequence of semigroups
in strong operator sense. However, in many cases, semigroups and
their associated closed  forms may live on different Hilbert spaces.
Fortunately, the Mosco convergence theory can be extended to cover
these cases of varying state spaces. Theorem \ref{t:equi} in the
Appendix,
which was obtained in \cite{K0} and \cite[Theorem 2.5]{K},
gives
one such extension. See \cite{KS} for another extension.

 In this
section, we establish the Mosco convergence of $(\sE^{(k)},
\sF^{(k)})$ in the sense of Definition \ref{D:Mosco}.
  For this, we  define the restriction operator $\pi_k:
L^2(M,m)\to L^2(V_k,m_k)$  and the extension operator $E_k:
L^2(V_k,m_k) \to L^2(M,m)$ as follows:
\begin{eqnarray*}
\pi_k f(x)&=& \frac 1{m_k(x)}\int_{U_k(x)}f(y)m(dy)\qquad
\hbox{for } f\in L^2(M,m)
\hbox{ and }  x\in V_k,\\
E_kg(z)&=& g(x)\qquad \hbox{for } g\in L^2(V_k,m_k) \ \hbox{ and } \
z\in  \mbox{Int}\,  U_k(x) \hbox{ with }  x\in V_k.
\end{eqnarray*}

Let $\langle \cdot, \cdot \rangle_k$ (resp. $\langle \cdot, \cdot
\rangle$) be the  inner product in Hilbert space $L^2(V_k,m_k)$
(resp. $L^2(M,m)$) and
 $\|\cdot\|_{k,p}$ (resp. $\|\cdot\|_p$)
be the $L^p$-norm of $L^p(V_k,m_k)$ (resp. $L^p(M,m)$).

\begin{lemma} \label{l:sspi1} \begin{description}
\item{\rm (i)} $\pi_k$ is a bounded operator from  $L^2(M,m)$  to  $L^2(V_k,m_k)$ with
$\sup_{k \ge 1} \|\pi_k\|  \le 1$, where $\|\pi_k\|$ is the
operator norm of $\pi_k$.
 Further, $\lim_{k \to \infty} \|\pi_k f
\|_{k,2} = \|f\|_2$ for every $f\in L^2(M; m)$.

\item{\rm (ii)} For each $f_k \in L^2(V_k,m_k)$, we have the following;
\begin{eqnarray}
 \pi_k E_k f_k &=& f_k
 \qquad m\text{-a.e.},
 \label{e:piE0}\\
 \langle \pi_k g, f_k \rangle_k &=& \langle g, E_k f_k \rangle \qquad\hbox{for every }
   g \in  L^2(M, m).   \label{e:piE1}
   \end{eqnarray}

\item{\rm (iii)} For every $f \in  L^2(V_k,m_k)$, $E_kf \in L^2(M, m)$ and
$\|E_kf \|_{2}^2 =  \|E_k(f^2) \|_{1} = \|f \|^2_{k,2}.
$

\item{\rm (iv)}
 For every $f\in L^2(M; m)$,
$E_k \pi_k f$ converges strongly to $f$ in $L^2 (M, m)$.

\item{\rm (v)} Suppose $f \in  C_c(M)$. Let $f_k:=f|_{V_k}\in L^2(V_k, m_k)$.
Then $E_k f_k$ converges strongly to $f$ in $L^2 (M, m)$.
\end{description}
\end{lemma}

\proof  (i) By the Cauchy-Schwarz inequality,
\begin{eqnarray}
\|\pi_kf\|^2_{k,2} & =& \sum_{x  \in V_k}m_{k}(x)
\left(\frac 1{m_{k}(x)}\int_{U_{k}(x)}f(y)m(dy)\right)^2 \label{e:fads}\\
&\le& \sum_{x  \in V_k}\frac {m_{k}(x)}{m_{k}(x)}
\int_{U_{k}(x)}f(y)^2m(dy)\, =\,\|f\|_2^2. \nn
\end{eqnarray}
 Moreover, by the uniform continuity, we
 easily see from \eqref{e:fads}
  that $\lim_{k \to \infty}\|\pi_k f \|^2_{k,2}=\|f\|^2_2$
for $f\in C_c(M)$.
As $C_c(M)$ is dense in $L^2(M; m)$ and
$\| \pi_k\|\leq 1$,
we have $\lim_{k \to \infty}\|\pi_k f \|^2_{k,2}=\|f\|^2_2$
for $f\in L^2(M; m)$.

(ii) \eqref{e:piE0} is clear from the definitions of  $\pi_k$ and
$E_k$.  The left hand side of \eqref{e:piE1} is
$$\sum_{x \in V_k}
\frac 1{m_{k}(x)}\int_{U_{k}(x)}g(y)m(dy)f_k(x)m_{k}(x).
$$
By Fubini's theorem, the above is equal to
\[ \int_{M}\sum_{x \in V_k}f_k(x)g(y)\1_{U_{k}(x)}(y)m(dy)\,=\,
\langle  E_k f_k,  g \rangle.
\]

(iii) Note that,
since $m(U_k(x)\cap U_k(y))=0$ for $x\ne y$,
 we have for $f \in  L^2(V_k,m_k)$
\begin{eqnarray*}
\|E_k f \|_2^2&=&   \int_{M}
\Big(\sum_{x \in V_k}f(x)\1_{U_{k}(x)}(y)\Big)^2 m(dy)\\
&=& \int_{M}\sum_{x\in V_k}f(x)^2\1_{U_{k}(x)}(y) m(dy)\,=\,
\|E_k(f^2) \|_1.
\end{eqnarray*}
Moreover, by Fubini's theorem,
$$
\int_{M}\sum_{x \in V_k}f(x)^2\1_{U_{k}(x)}(y) m(dy) \,= \,\sum_{x
\in V_k} f(x)^2 m_k(x) =\|f \|^2_{k,2}.
$$

(iv) First assume that $f\in C_c(M)$. Let $K:=\{x\in M: \rho (x, {\rm supp}[f])\leq 1\}$.
 By the Cauchy-Schwarz inequality, for sufficiently large $k\geq 1$,
\begin{eqnarray*}
\|E_k f-f \|_2^2
&=&\int_{K} |E_kf_k(x)-f(x)|^2 \, m(dx)  \\
&\le& \sum_{z\in V_k\cap K}
\int_{U_k(z)} \left(
 \frac{1}{m_k(z)} \int_{U_k(z)} (f(y) - f(x)) m(dy)  \right)^2
m(dx)\\
 &\leq & \sum_{z\in V_k\cap K}  \frac{1}{m_k(z)} \int_{U_k(z)\times
 U_k(z)} (f(y) - f(x))^2 m(dy) m(dx) ,
\end{eqnarray*}
which, by the uniform continuity of $f\in C_c(M)$,  tends to zero as
$k\to \infty$.  That is, for $f\in C_c(M)$, $E_k \pi_k f$ converges
strongly in $L^2(M; m)$ to $f$. Since by (i) and (iii),
$$\|E_k \pi_k f\|_2 = \| \pi_k f\|_{k, 2} \leq \| f\|_2
\qquad \hbox{for } f\in L^2(E; m)
$$
and that $C_c(M)$ is dense in $L^2(E; m)$, we conclude that for every $f\in L^2(E; m)$,
$E_k \pi_k f$ converges
strongly in $L^2(M; m)$ to $f$.

(v)  Let $K:=\{x\in M: \rho (x, {\rm supp}[f])\leq 1\}$. Then for $k$ sufficiently large,
$$
\int_{M} |E_kf_k(x)-f(x)|^2 m(dx)
=
 \int_{K} |E_kf_k(x)-f(x)|^2 m(dx), $$
 which goes to zero  by the uniform continuity of $f$.
\qed

\subsection{Mosco convergence}

Next we assume the following:

\medskip

\noindent {\bf (A2)}. {\it For $m$-a.e. $x\in M$, $j(x,\cdot)$ is a  positive
 measure on $M\setminus \{x\}$ such that the following holds:
\begin{description}
\item{\rm (i)}
  For any $\eps>0$,
 $x\mapsto j(x, M\setminus B(x,\eps))$
is locally integrable with respect to $m$.

\item{\rm (ii)}
 For any non-negative Borel measurable functions
$u,v$,
\[\int_Mu(x)(jv)(x)m(dx)=\int_M(ju)(x)v(x)m(dx)\,(\le \infty).\]
Here $ju(x):=\int_{M\setminus \{x\}} u(y) j(x, dy)$.

\item{\rm (iii)}
For any compact set $K$,
\begin{equation}\label{e:A1}
\sup_{x\in K} \int_{M} (\rho(x,y) \wedge 1)^2 j(x, dy)
  <\infty .
\end{equation}
\end{description}
} \noindent
Let  ${\wh d}$ denote the diagonal set in $M\times M$.
The
kernel $j$ then determines a positive Radon measure $J(dx,dy)$ on
$M\times M\setminus \wh d$  by
\[
 \int_{M\times M\setminus \wh d} f(x,y)J(dx,dy)= \int_M \left( \int_M
 f(x,y)j(x,dy) \right) m(dx) \quad \hbox{for }  f\in C_c(M\times M \setminus \wh d).
\]
Define a bilinear form $(\sE,\sF)$ on $L^2(M; m)$ as follows:
\begin{eqnarray}
\sF&:=& \left\{ u\in L^2(M, m): \,\int_{M\times M \setminus \wh d}{(u(x)-u(y))^2}
J(dx,dy)<\infty \right\}, \label{e:Df} \\
\sE (u, v)&:=& \frac1{2} \int_{M\times M \setminus \wh d}
{(u(x)-u(y))(v(x)-v(y))} {J(dx,dy)} \qquad \hbox{for } u, v\in \sF.
\nn
\end{eqnarray}

\begin{lemma}\label{lind}
Under the condition
{\bf (A2)}, $\Lip_{c} (M)\subset \sF$.
\end{lemma}

\proof Let $u\in \Lip_c(M)$. Clearly it is $L^2(M, m)$-integrable.
Denote by $\Lambda$ the  Lipschitz constant of $u$ and $K:={\rm
supp}[u]$. Then by the symmetry of
$j(x,dy)$,
\begin{eqnarray*}
\sE (u, u) &\leq & \int_K \left( \int_{M\setminus \{x\}}
(u(x)-u(y))^2 j(x,  dy)\right) m(dx) \\
&\leq& \int_K \left( \int_{M} \left( \Lambda^2 \rho (x, y)^2 {\bf
1}_{\{\rho (x, y) \leq 1\}}+ 4 \| u\|_\infty^2 {\bf 1}_{\{\rho (x,
y) >1 \}} \right) j(x, dy) \right) m(dx) \\
&\leq & C \, m(K) \, \sup_{x\in K} \int_M  \left( \rho (x, y)^2
\wedge 1 \right)  j(x, dy) ,
\end{eqnarray*}
which is finite by condition \eqref{e:A1}. This proves that $u\in
\sF$. \qed

Lemma \ref{lind} in particular implies that $\sF$ is a dense linear
subspace of $L^2(M; m)$. It is easy to check by using Fatou's lemma
that $(\sE,\sF)$ is Dirichlet form on $L^2(M; m)$ (cf. \cite[Example
1.2.4]{FOT}). We further assume that

\medskip

\noindent{\bf (A3)}. {\it $\Lip_c(M)$ is dense in $(\sF,
\sE(\cdot,\cdot)+\|\cdot\|_2^2)$. }
\medskip

Under conditions {\bf (A2)} and {\bf (A3)},  $(\sE,\sF)$ is a
regular Dirichlet form on $L^2(M;m)$. Let $X=\{X_t, t\geq 0, \P_x,
x\in M\}$ be its associated symmetric Hunt process on $M$. If
\begin{equation}\label{e:cons2}
\sup_{x\in M} \int_M (\rho (x, y)^2\wedge 1)j(x, dy)<\infty,
\end{equation}
 then we have by   \eqref{e:2.2} and \cite[Theorem
3.1]{MU} that the process $X$ is conservative; that is, $X$ has
infinite lifetime $\bP_x$-a.s. for $\sE$-q.e. $x\in M$.

In the following, we sometimes extend $\{\sC^{(k)}(x,y): x,y\in V_k\}$ to
$\{\sC^{(k)}(z,w): z\in \mbox{Int}\,U_{k}(x),\, w\in \mbox{Int}\,U_{k}(y),\, x,y\in V_k\}$ by taking
$\sC^{(k)}(z,w)=\sC^{(k)}(x,y)$ for $z\in \mbox{Int}\,U_{k}(x)$ and $w\in \mbox{Int}\,U_{k}(y)$.
We recall that we have fixed some $x_0\in M$ and $B_r=B(x_0, r)$.
We will use the following definition for the remainder of this paper:
Define for
function $f: M\to \bR$,
\begin{equation}\label{e:jdeltak}
\overline{\sE}^{(k)}_{j, \delta} (f,f):=\frac1{2} \int\int_{\{ (z,w) \in B_j
\times B_j: \,  \rho(z,w) > \delta\}}(f(w)-f(z))^2 \sC^{(k)}(w,z)
m(dw)m(dz),
\end{equation}
\begin{equation}\label{e:jdelta}
\sE_{j, \delta} ( f,  f) :=
\frac12 \int\int_{\{(z,w) \in B_j \times B_j: \,  \rho(z,w) >
\delta \}}(f(w)-f(z))^2J(dw,dz)
\end{equation}
and
$$
\overline{\sE}^{(k)} (f,f):=\frac1{2} \int_{M\times M}
(f(w)-f(z))^2  \sC^{(k)}(w,z)
m(dw)m(dz).
$$
Note that for function $f$ on $V_k$,
$
(E_kf(z) -E_kf(w))^2 = (f(x)- f(y))^2
$ where
$x,y \in V_k$ with $z\in U_k(x), w\in U_k(y)$.
Thus
\bee\label{e:dvd1}
 \overline{\sE}^{(k)}(E_k u,\,  E_k u)
 \, =\, \sE^{(k)} (u, u), \quad \text{ for all } u \in \sF^{(k)}
\eee

Our final assumption in this subsection is the following.

\medskip

{\it \noindent {\bf (A4)}.
{\rm (i)} For any compact subset $K\subset M$,
 \begin{eqnarray} \label{e:A4.1}
 && \lim_{\eta\to 0} \limsup_{k\to\infty}
 \int\int_{\{(x, y)\in K\times K: \rho(x, y)\leq  \eta\}} \rho(x, y)^2 \sC^{(k)}(x,y)
   m(dx)m(dy)=0,
\end{eqnarray}
 \begin{eqnarray} \label{e:A4.2}
   \lim_{j\to \infty} \limsup_{k\to\infty} \int_K
 \int_{B_j^c} \sC^{(k)} (x, y) m(dx)m(dy)
 =0.
\end{eqnarray}

{\rm (ii)} For every $\eps>0$, there exists $N >0$ such that
for every $k \ge i\ge N$ and
$f\in L^2(V_i; m_i)$,
$$ \sE^{(k)} ( \pi_k E_i f, \, \pi_k E_i f)^{1/2} \leq \sE^{(i)} (f, f)^{1/2}+\eps.
$$

{\rm (iii)} For any sufficiently small $\delta>0$ and large $j\in \bN$,
  \begin{eqnarray}\label{e:A4.3}
   \lim_{k\to \infty}
\overline{\sE}^{(k)}_{j, \delta} (f,f) = {\sE}_{j, \delta} (f,f)
  \qquad \hbox{for every } f \in \Lip_c(M).
  \end{eqnarray} }

\begin{remark}\label{R:4.2} \rm
 It follows from \eqref{e:A1} of {\bf (A2)}
 that for every compact subset $K\subset M$
\begin{equation}\label{e:4.7}
 \lim_{\eta \to 0}  \int_{\{(x, y)\in K\times M: \rho(x, y)\leq \eta\}} \rho (x,
  y)^2 j(x,dy) m(dx)=0.
\end{equation}
 \qed
\end{remark}

\begin{lemma}\label{l:fcon1}
Suppose the conditions {\bf
(A2)}, {\bf
(A3)} and {\bf (A4)} {\rm (i)(iii)} hold, then
for every $f \in \Lip_c(M)$, $\lim_{k\to \infty} \sE^{(k)}(\pi_k f,
\pi_k f)=\sE(f,f)$.
\end{lemma}

\proof
First, note that by (\ref{e:dvd1}),
$\sE^{(k)}(\pi_k f,\pi_k f)=\overline{\sE}^{(k)}(E_k \pi_k f,\,  E_k \pi_k f)$.

 Fix $f \in \Lip_c (M)$ and let $K$ be the support of $f$, $K_1:=\{x \in M: \rho(x,K)
\le 1  \}$ and $M_f:=\sup_{x\in M} |f(x)|$. Then, by (\ref{e:A4.2}) and
the symmetry of
$\sC^{(k)}$
for each $\eps>0$, there exists $j_0$ such that the following
holds for $j\ge j_0$,
\begin{eqnarray*}
&&\limsup_{k\to\infty}\frac 12 \int\int_{(B_j\times B_j)^c} (E_k \pi_k f(x)-E_k \pi_k f(y))^2\sC^{(k)} (x, y) m(dx)m(dy)\\
 &\le & (2M_f)^2
 \limsup_{k\to\infty}\int_{K}
 \int_{B_j^c}\sC^{(k)} (x, y) m(dx)m(dy)<\eps .
 \end{eqnarray*}
Similarly, using (\ref{e:A1}) and choosing $j_0$ larger if necessary,
we have
\[\frac 12 \int\int_{(B_j\times B_j)^c} (f(x)-f(y))^2J(dx,dy)<\eps.\]
Since $f \in \Lip_c(M)$ is Lipschitz continuous,
using {\bf (AG.2)}, {\bf (AG.3)}, (\ref{e:A4.1}) and
\eqref{e:4.7} and
arguing similarly,  we have
$$
\limsup_{k\to\infty}\frac 12 \int
 \int_{\{(x,y)\in K_1\times K_1: \rho(x,y) \le  \delta\}} (E_k \pi_k f(x)-E_k \pi_k f(y))^2\sC^{(k)} (x, y) m(dx)m(dy)<\eps$$
 and
 $$\frac 12 \int
 \int_{\{(x,y)\in K_1\times K_1: \rho(x,y) \le  \delta\}}(f(x)-f(y))^2J(dx,dy)<\eps
 $$
for all $\delta \in (0,1)$. Thus, it is enough to show
the following for any sufficiently small $\delta$ and large $j$:
\bee\label{e:sjoebf}
\lim_{k\to\infty}\overline{\sE}^{(k)}_{j,\delta}(E_k \pi_k f,\,  E_k \pi_k f)
=\sE_{j, \delta}(f,f).
\eee

  By the symmetry of $\sE^{(k)}_{j,\delta}$ and Lemma \ref{l:sspi1} (iv),
\begin{eqnarray*}
&& \lim_{k\to \infty} \left|\overline \sE^{(k)}_{j,\delta}(E_k \pi_k f, E_k \pi_k f)^{1/2}-
\overline{\sE}^{(k)}_{j,\delta}( f, f)^{1/2}\right| \nonumber \\
&\leq & \lim_{k\to \infty} \overline \sE^{(k)}_{j,\delta}(f-E_k \pi_k f,\, f-E_k \pi_k f)^{1/2} \nonumber \\
&=& \lim_{k\to \infty}
\left( \frac12
 \int_{B_j\times B_j} \big((f-E_k\pi_kf)(x)-(f-E_k\pi_kf)(y)\big)^2
 \sC^{(k)}(x,y) \1_{\{\rho (x, y)>\delta\}}m(dx)m(dy) \right)^{1/2} \nonumber \\
 &\leq & \lim_{k\to \infty}
 \left(
 \int_{B_j} (f(x)-E_k \pi_k f(x))^2 \left( \int_{B_j}
    \sC^{(k)}(x,y) \1_{\{\rho (x, y)>\delta\}}m(dy)\right)m(dx)\right)^{1/2} \nonumber \\
    &\leq & \lim_{k\to \infty}
    c(j, \delta)
     \, \| f-E_k \pi_k f\|_2 =0.
\end{eqnarray*}
 Hence we have
  \begin{equation}\label{e:4.15}
 \lim_{k\to \infty} \overline\sE^{(k)}_{j,\delta}(E_k \pi_k f, E_k \pi_k f)=\lim_{k\to \infty}
\overline{\sE}^{(k)}_{j,\delta}( f, f).
\end{equation}
On the other hand,
by \eqref {e:A4.3},
  $\lim_{k\to \infty} \overline \sE^{(k)}_{j,\delta}(f,f) =\sE_{j, \delta}(f, f)$.
 This completes the proof of the Lemma.
\qed

The following lemma is needed in establishing the Mosco convergence of
$(\sE^{(k)}, \sF^{(k)})$ to $(\sE, \sF)$. It is formulated in a general setting.

\begin{lemma}\label{l:scmc1}
Suppose $(\sH_k,  \langle \,\cdot\,, \,\cdot\, \rangle_k )$ and
$(\sH,  \langle \,\cdot\,, \,\cdot\, \rangle)$ are Hilbert spaces
with corresponding norms $\|\cdot\|_k$ and $\| \cdot \| $ respectively.
Suppose that for each $k\geq 1$, there is a bounded linear operator
$E_k: \sH_k\to \sH$ so that its adjoint operator $E_k^*=:\pi_k$  is its left
inverse  satisfying the conditions \eqref{eqn:AAA}--\eqref{eqn:A} in the Appendix.
Let $a^{(k)}$ be a
symmetric bilinear form on $\sH_k$ and $a$ be a symmetric bilinear
form on $\sH$. Then,
$a^{(k)}$ is Mosco convergent to $a$ in the generalized sense of
Definition \ref{D:Mosco}
 if Definition \ref{D:Mosco}(i) holds and in addition the
following hold:
 \begin{description}
\item{({1})}
There exists a set $\sD \subset \sH$ which is dense in $(\sD[a], a+ \|\cdot\|^2)$.
\item{({2})} $\pi_k (\phi) \in \sD[a^{(k)}]$ for every $\phi \in \sD$.
\item{({3})} For every $\phi \in \sD$,
$$
 \limsup_{k \to \infty} \, a^{(k)} (\pi_k \phi, \pi_k\phi) = a(\phi,\phi).
$$
\end{description}
\end{lemma}
\proof
Note that,  since $\| u_k\|_k = \|E_k u_k\|$, the strong convergence of
$ E_k u_k$ to  $u$ in $\sH$ is equivalent to that
$\| u_k\|_k \to \|u\|$ and the weak convergence of
$ E_k u_k$ to  $u$ in $\sH$.
Thus
the proof of this lemma is the same as the one of \cite[Lemma 2.8]{Ko}.
\qed

\begin{theorem}\label{Mosco1}
Suppose
the conditions \eqref{e:A1a} of {\bf (A1)} and {\bf (A2)}--{\bf (A4)} hold, then
$(\sE^{(k)}, \sF^{(k)})$ is Mosco convergent to $(\sE, \sF)$ in the generalized sense of
Definition \ref{D:Mosco}.
\end{theorem}

\proof Take $\sD=\Lip_c(M)$ in Lemma \ref{l:scmc1}. Then,  By our
assumption {\bf (A3)} and  Lemmas \ref{lind}, \ref{l:fcon1},
\ref{l:scmc1}, we only need to check   condition (i) in Definition
\ref{D:Mosco}.

It is enough to consider sequences $\{u_k\}_{k \ge 1} \subset
L^2(V_k, m_k)$ such that $E_ku_k$ converges weakly  to $u \in
L^2(M,m)$ and $\liminf_{k\to \infty} \sE^{(k)} ( u_k,  u_k)  <
\infty$.
Taking a subsequence if necessary, we may and do assume
that  $\lim_{k \to \infty} \sE^{(k)} (u_{k}, u_{k} )$ exists and is
finite, and that
  \bee\label{e:fsup1} \sup_{k \ge 1} \, \left(
\sE^{(k)} ( u_k,  u_k) + \sum_{x\in V_k} u_k(x)^2 m_k(x)\right) < \,
\infty .
  \eee
   So in particular, $ u_k \in \sF^{(k)}$ for every $k\geq 1$.
By uniform boundedness principle,
  $\{E_k u_k; \, k \ge 1\}$ is a bounded sequence on $L^2(M; m)$.

By the Banach-Saks theorem, taking a subsequence if necessary,
$v_k:=\frac{1}{k}\sum^{k}_{i=1}E_i u_{i} $ converges to some
$v_\infty$ in  $L^2(M; m)$. Since $E_k u_k$ converges weakly to $u$
in $L^2(M; m)$,  $v_\infty $ must be $u$ $m$-a.e. on $M$.

Fix an integer $j\geq 1$ and $\delta>0$. For $\eps >0$, let $f\in \Lip_c(M)$
such that $\| u-f\|_2 \leq \eps /
\sqrt{2a_{j, \delta}} $,
where
$$ a_{j, \delta} = \max\left\{\sup_{k \geq k_0} \sup_{x \in B_j } \int_{M}  \sC^{(k)}(x,y) \, \1_{\{\rho (x, y)>\delta\}} \, m (dy),  \
   \sup_{z \in B_j }  \int_{M}  \1_{\{\rho (x, y)>\delta\}} \, j(z,dw)\,  \right\},
 $$
which is finite by \eqref{e:A1a} of {\bf (A1)} and {\bf (A2)}(iii).
Observe that by {\bf (A4)}(iii)
\begin{eqnarray*}
  && \limsup_{k\to \infty} \Big| \overline \sE^{(k)}_{j, \delta} (v_k, v_k)^{1/2} -  \sE_{j, \delta} (f, f)^{1/2} \Big| \\
 &\leq &  \limsup_{k\to \infty} \Big| \overline \sE^{(k)}_{j, \delta} (v_k, v_k)^{1/2} -  \overline \sE^{(k)}_{j, \delta} (f, f)^{1/2} \Big| \\
 &\leq & \limsup_{k\to \infty}  \overline \sE^{(k)}_{j, \delta} (v_k-f, v_k-f)^{1/2} \\
 &\leq &
   \limsup_{k\to \infty} \left( 2
 \int_{B_j} (v_k(x)-f(x))^2 \left(\int_M \sC^{(k)}(x, y)
  \1_{\{\rho (x, y)>\delta\}} m(dy) \right) m(dx) \right)^{1/2} \\
  &\leq &  \limsup_{k\to \infty}
    \sqrt{2a_{j, \delta}}
  \, \| v_k -f\|_2
   =
   \sqrt{2a_{j, \delta}}
    \, \| u-f\|_2 <\eps.
\end{eqnarray*}
Similarly, we have
$$ \Big| \sE_{j, \delta} (f, f)^{1/2} -\sE_{j, \delta} (u, u)^{1/2} \Big|
\leq \sE_{j, \delta} (f-u, f-u)^{1/2}
\leq
   \sqrt{2a_{j, \delta}}\, \|f-u\|_2<\eps.
$$
Thus we have
\begin{equation}\label{e:ddf3}
   \liminf_{k\to \infty}  \overline \sE^{(k)}_{j, \delta} (v_k, v_k)^{1/2}  \geq   \sE_{j, \delta} (f, f)^{1/2}-\eps  \geq   \sE_{j, \delta} (u, u)^{1/2}-2\eps.
\end{equation}

Observe that for $k_0\le n \le k$,
\begin{eqnarray*}
   && \overline \sE^{(k)}_{j, \delta} (v_n, v_n)^{1/2} -  \sE_{j, \delta} (f, f)^{1/2}  \\
 &\leq &  \Big| \overline \sE^{(k)}_{j, \delta} (v_n, v_n)^{1/2} -  \overline \sE^{(k)}_{j, \delta} (f, f)^{1/2} \Big| \\
 &\leq &   \overline \sE^{(k)}_{j, \delta} (v_n-f, v_n-f)^{1/2} \\
 &\leq &
 \left(2\int_{B_j} (v_n(x)-f(x))^2 \left(\int_M \sC^{(k)}(x, y)
  \1_{\{\rho (x, y)>\delta\}} m(dy) \right) m(dx) \right)^{1/2} \\
  &\leq &
   \sqrt{2a_{j, \delta}} \, \| v_n -f\|_2.
\end{eqnarray*}
Thus
$$
\lim_{n \to  \infty} \sup_{k \ge n} \overline \sE^{(k)}_{j, \delta} (v_n, v_n)^{1/2} \le
\sE_{j, \delta} (f, f)^{1/2} +
   \sqrt{2a_{j, \delta}}\, \| u-f\|_2 \le
   \sE_{j, \delta} (f, f)^{1/2} +\eps < \infty.
$$

By condition {\bf (A4)}(ii) and the above, there exists $N>0$ such that
for every $k \ge i\ge N$,
\begin{equation}\label{e:ddf1}
  \sE^{(k)} ( \pi_k E_i u_i, \, \pi_k E_i u_i)^{1/2} \leq \sE^{(i)} (u_i, u_i)^{1/2}+\eps ;
\end{equation}
and
\begin{equation}\label{e:ddf2}
\sup_{m \ge N} \overline \sE^{(m)}_{j, \delta} (v_N, v_N)^{1/2}  < \infty.
 \end{equation}
Since, for $k >N$
\begin{eqnarray*}
  \overline \sE^{(k)}_{j, \delta} (v_k, v_k)^{1/2}
&=&
 \overline \sE^{(k)}_{j, \delta} \Big(\frac{1}{k}\sum^k_{i=1}E_i u_i , \,
 \frac{1}{k}\sum^k_{i=1}E_i u_i \Big)^{1/2}\\
 &=&
 \overline \sE^{(k)}_{j, \delta} \Big(\frac{1}{k}\sum^N_{i=1}E_i u_i +\frac{1}{k}\sum^k_{i=N+1}E_i u_i , \,
\frac{1}{k}\sum^N_{i=1}E_i u_i +\frac{1}{k}\sum^k_{i=N+1}E_i u_i  \Big)^{1/2}\\
 &\le& \frac{N}{k} \overline \sE^{(k)}_{j, \delta} (v_N, v_N)^{1/2}+ \frac1{k}\sum^k_{i=N+1}
 \overline{\sE}^{(k)}_{j, \delta} \big(E_i u_i , \, E_i u_i\big)^{1/2}\\
  &\le& \frac{N}{k} \left(\sup_{m \ge N}\overline \sE^{(m)}_{j, \delta} (v_N, v_N)^{1/2} \right)+ \frac1{k}\sum^k_{i=N+1}
{\sE}^{(k)} \big(\pi_k E_i u_i , \, \pi_k E_i u_i\big)^{1/2}
 \end{eqnarray*}
by \eqref{e:ddf1}--\eqref{e:ddf2},
\begin{eqnarray*}
 \liminf_{k\to \infty}  \overline \sE^{(k)}_{j, \delta} (v_k, v_k)^{1/2}
&\le& \liminf_{k\to \infty}  \frac1{k}\left(\sum^k_{i=N+1}
{\sE}^{(i)} \big(u_i , \, u_i\big)^{1/2}\right)+\eps \\
&\le& \lim_{k \to \infty}
 \sE^{(k)}( u_{k} , \,   u_{k})^{1/2}+\eps.
 \end{eqnarray*}
 Now from \eqref{e:ddf3}, we have
 $$
 \sE_{j, \delta} (u, u)^{1/2} \le
  \lim_{k \to \infty}
 \sE^{(k)}( u_{k} , \,   u_{k})^{1/2}+3\eps.
 $$
 Since $\eps>0$ is arbitrary, we have
 $$
 \sE_{j, \delta} (u, u) \le
  \lim_{k \to \infty}
 \sE^{(k)}( u_{k} , \,   u_{k}).
 $$
By first letting $j\to \infty$ and then $\delta \to 0$, one has $\lim_{k \to \infty}
 \sE^{(k)}( u_{k} , \,   u_{k}) \geq \sE (u, u)$, which completes the
  proof of the theorem. \qed

\subsection{Mosco convergence under alternative setup}

We first give an alternative assumption and give Mosco convergence.
We do {\it not} assume {\bf (A1)} in this subsection.
For $u\in L^2(B_j, m)$, define
\begin{eqnarray*}
\overline \sL^{(k)}_{j,\delta}u(x)&=&
\int_{B_j}(u(y)-u(x))
\sC^{(k)}(x,y)1_{\{\rho(x,y)>\delta\}}m(dy)~~~~\forall x\in B_j,\\
\sL_{j,\delta}u(x)&=&
\int_{B_j}(u(y)-u(x))1_{\{\rho(x,y)>\delta\}}j(x,dy)~~~~\forall x\in B_j.
\end{eqnarray*}
Then $\overline\sE^{(k)}_{j,\delta}(u,v)=-(u,\overline \sL^{(k)}_{j,\delta}v)_{2,B_j}$ and
$\sE_{j, \delta}(u,v)=-(u,\sL_{j,\delta}v)_{2,B_j}$
where $(u,v)_{2,B_j}
=\int_{B_j} u(x)v(x)m(dx)$
and $\overline\sE^{(k)}_{j,\delta}(u,v)$ and $\sE_{j, \delta}(u,v)$ are defined in
\eqref{e:jdeltak} and  \eqref{e:jdelta} respectively.
\\

In this subsection,
we assume {\bf (A2)}, {\bf (A3)${}^*$}
and {\bf (A4)${}^*$} below:

{\it \noindent {\bf (A3)${}^*$.}\,{\rm (i)} Same as {\bf (A3)} in Section \ref{sec4}.

{\rm (ii)}
$\sL_{j,\delta} f$ is continuous for all $f\in \Lip_c(M).$}

{\it \noindent {\bf (A4)${}^*$.}\,{\rm (i)} Same as {\bf (A4)}{\rm (i)} in Section \ref{sec4}.

{\rm (ii)}
For any sufficiently small $\delta>0$ and large $j\in \bN$,
\[\lim_{k\to\infty}\int_{B_j}(\overline \sL^{(k)}_{j,\delta}f(x))^2m(dx)
=\int_{B_j}(\sL_{j,\delta} f(x))^2m(dx),~~~\forall f\in \Lip_c(M).\]
\indent {\rm (iii)}
For any sufficiently small $\delta>0$ and large $j\in \bN$,
\begin{eqnarray*}
  \lim_{k\to \infty}
\overline{\sE}^{(k)}_{j, \delta} (f,f) = {\sE}_{j, \delta} (f,f)
   \qquad \hbox{for every } f\in C_b(B_j).
  \end{eqnarray*} }
In other words, we put additional assumption {\bf (A3)${}^*$}\,{\rm (ii)},
 strengthen {\bf (A4)}\,(iii),
 and
replace {\bf (A4)}\,(ii) in Section \ref{sec4} by {\bf (A4)${}^*$}\,(ii).
Note that, by the polarization identity,  {\bf (A4)${}^*$} (iii) is equivalent to
\begin{equation} \label{e:polar}
 \lim_{k\to \infty}
\overline{\sE}^{(k)}_{j, \delta} (f,g) = {\sE}_{j, \delta} (f,g)
   \qquad \hbox{for every } f, g\in C_b(B_j).
  \end{equation}

Let
\[\sup_{x\in B_j}\int_{B_j}1_{\{\rho(x,y)>\delta\}}j(x,dy)=:K_{j,\delta},\]
which is finite due to {\bf (A2)}.
Also, let $\|\cdot\|_{2,B_j}$ be the $L^2$-norm on $B_j$.
We then have the following basic estimates.
\begin{lemma}\label{l:bas_lf}
The following holds for any $\delta>0$ and $j\in \bN$.\\
{\rm (i)} ~$\sE_{j, \delta} (u,u)\le K_{j,\delta}\|u\|_{2,B_j}^2$ for all $u\in L^2(B_j, m)$.
Especially, $\sE_{j, \delta} (u,u)<\infty$ for all $u\in L^2(B_j, m)$.\\
{\rm (ii)}~$\|\sL_{j,\delta} u\|_{2,B_j}^2\le  K_{j,\delta} \sE_{j, \delta} (u,u)$ for all $u\in L^2(B_j, m)$.\\
{\rm (iii)} $\lim_{k\to\infty}\|(\sL_{j,\delta} -\overline \sL^{(k)}_{j,\delta})f\|_{2,B_j}=0$ for all
$f\in \Lip_c(M)$.
\end{lemma}
\proof (i) For $u\in L^2(B_j, m)$, we have
\begin{eqnarray*}
\sE_{j, \delta} (u,u)&=&\frac 12 \int_{B_j}\int_{B_j}
(u(x)-u(y))^2j(x,y)1_{\{\rho(x,y)>\delta\}}dxdy\\
&\le & \|u\|_{2,B_j}^2\sup_{x\in B_j}\int_{B_j}j(x,y)1_{\{\rho(x,y)>\delta\}}dy\le K_{j,\delta}\|u\|_{2,B_j}^2.
\end{eqnarray*}
(ii) As in (i), $\sE_{j, \delta} (u,u)<\infty$ for $u\in L^2(B_j, m)$. So,
using the Cauchy-Schwarz inequality, we have
\begin{eqnarray*}
\|\sL_{j,\delta} u\|_{2,B_j}^2&=& \int_{B_j}\Big(\int_{B_j}(u(y)-u(x))1_{\{\rho(x,y)>\delta\}}j(x,dy)\Big)^2m(dx)\nonumber\\
&\le &  \int_{B_j}
\Big(\int_{B_j}(u(x)-u(y))^21_{\{\rho(x,y)>\delta\}}j(x,dy)\cdot
\int_{B_j}1_{\{\rho(x,y)>\delta\}}j(x,dy)\Big)m(dx)~~~\\
&\le & K_{j,\delta} \sE_{j, \delta} (u,u).
\end{eqnarray*}
(iii) Using
{\bf (A3)${}^*$}(ii) and
{\bf (A4)${}^*$}(ii)(iii) (and \eqref{e:polar}), we have
\[\|(\sL_{j,\delta} -\overline \sL_{j,\delta}^{(k)})f\|_{2,B_j}^2=\|\sL_{j,\delta} f\|_{2,B_j}^2+\|\overline \sL_{j,\delta}^{(k)}f\|_{2,B_j}^2
-2\overline\sE^{(k)}_{j,\delta}(\sL_{j,\delta} f, f)\to 0.\]
\qed

We now prove the Mosco convergence that corresponds to Theorem \ref{Mosco1}.
Recall that we do not assume {\bf (A1)} in this subsection.
\begin{theorem}\label{Mosco1--0}
$(\sE^{(k)}, \sF^{(k)})$ is Mosco convergent to $(\sE, \sF)$ in the generalized sense of
Definition \ref{D:Mosco}.
\end{theorem}
\proof
Since
Lemma \ref{l:fcon1} works in this setting,
as before, we only need to check  condition (i) in Definition
\ref{D:Mosco}. Also, as in the proof of Theorem \ref{Mosco1},
we may assume $\{E_k u_k; \, k \ge 1\}$ is a bounded sequence on $L^2(M; m)$
that converges weakly  to $u \in L^2(M,m)$,
$\lim_{k\to \infty} \sE^{(k)} ( u_k,  u_k)  < \infty$, and (\ref{e:fsup1}) holds

In the following, we simply write $(\cdot,\cdot)$, $\|\cdot\|_2$ for
inner product and $L^2$-norm on $B_j$.
Fix $j$ large and $\delta>0$ small then take positive $\eps  <\sE_{j, \delta}(u,u)$. For $u\in L^2$ which is the weak limit of $E_ku_k$,
take $f\in \Lip_c(M)$ so that
$\sE_{j, \delta}(u-f,u-f)+\|u-f\|_2^2<\eps$ (note that
by Lemma \ref{l:bas_lf}(i), it is enough to take $\|u-f\|_2^2$ small).
First, note that
\bee\label{e:fstawe}
\lim_{k\to\infty}(E_ku_k,(\sL_{j,\delta} -\overline \sL_{j,\delta}^{(k)})f)=0,
\eee
where $u_k$, $u$ and $f$ are as above. Indeed, using Lemma \ref{l:bas_lf}(iii),
\[|(E_ku_k,(\sL_{j,\delta} -\overline \sL_{j,\delta}^{(k)})f)|\le \|E_ku_k\|_{2}
\|(\sL_{j,\delta} -\overline \sL_{j,\delta}^{(k)})f\|_{2}\le
\left(\sup_k \|E_ku_k\|_{2} \right)
\|(\sL_{j,\delta}
-\overline \sL_{j,\delta}^{(k)})f\|_{2}\to 0.\]

Now
\begin{eqnarray*}
|\overline\sE^{(k)}_{j,\delta}(E_ku_k,f)-\sE_{j, \delta}(f,f)|&=&
|(f,\sL_{j,\delta} f)-(E_ku_k, \overline \sL_{j,\delta}^{(k)}f)|\\
&\le&|(E_ku_k, (\sL_{j,\delta}-\overline \sL_{j,\delta}^{(k)})f)|+
|(E_ku_k-u, \sL_{j,\delta} f)|+|(u-f, \sL_{j,\delta} f)|.
\end{eqnarray*}
Using (\ref{e:fstawe}), the first term of the last line
goes to zero and since $\{E_ku_k\}$ converges weakly to $u$, the second term
goes to zero as $k\to\infty$
(note that $\sL_{j,\delta} f\in L^2$
due to Lemma
\ref{l:bas_lf}(ii)).
\ref{l:bas_lf}(i)(ii)).
Further,
there exists a
$C=C(j, \delta,u)>0$
such that
\begin{eqnarray*}
|(u-f, \sL_{j,\delta} f)|&\le& \|u-f\|_2\|\sL_{j,\delta} f\|_2\le \|u-f\|_2(\|\sL_{j,\delta} (u-f)\|_2+\|\sL_{j,\delta} u\|_2)\\
&\le &  \|u-f\|_2(K_{j,\delta}\|u-f\|_2+\|\sL_{j,\delta} u\|_2)\le C\eps^{1/2},
\end{eqnarray*}
where Lemma \ref{l:bas_lf}(i),(ii) are used in the third inequality.

Thus,
using Cauchy-Schwarz inequality,
we have
\begin{eqnarray*}
\sE_{j, \delta}(f,f)&\le&\limsup_{k\to\infty}|\overline\sE^{(k)}_{j,\delta}(E_ku_k,f)|+
C\eps^{1/2}\\
&\le& \lim_{k\to\infty}\Big(\overline\sE^{(k)}_{j,\delta}(E_ku_k,E_ku_k)^{1/2}
\sE^{(k)}_{j,\delta}(f,f)^{1/2}\Big)+C\eps^{1/2}
\\
&= & \lim_{k\to\infty}\overline\sE^{(k)}_{j,\delta}(E_ku_k,E_ku_k)^{1/2}
\overline\sE_{j, \delta}(f,f)^{1/2}+C\eps^{1/2}
\end{eqnarray*}
where the last equality is due to {\bf (A4)${}^*$}\,(iii).
Since $\eps <\sE_{j, \delta}(u,u)$,
by a rearrangement,  we obtain
\begin{eqnarray*}
\sE_{j, \delta}(u,u)^{1/2}&\le& \sE_{j, \delta}(f,f)^{1/2}+\eps^{1/2}
\le \lim_{k\to\infty}\overline\sE^{(k)}_{j,\delta}(E_ku_k,E_ku_k)^{1/2}+
C\frac{\eps^{1/2}}{\sE_{j, \delta}(f,f)^{1/2}}+\eps^{1/2}\\
&\le &\lim_{k\to\infty}\overline\sE^{(k)}(E_ku_k,E_ku_k)^{1/2}
+C\frac{\eps^{1/2}}{\sE_{j, \delta}(u,u)^{1/2}-\eps^{1/2}}+\eps^{1/2}.
\end{eqnarray*}
Taking $\eps\to 0$ and then $j\to\infty$ and $\delta\to 0$, we obtain the desired inequality.
\qed

\begin{remark} \label{r:n}
{\rm The assumption
{\bf (A3)}${}^*$\,(ii) is used only in the proof of Lemma \ref{l:bas_lf} (iii). Thus if we strengthen {\bf (A4)}${}^*$\,(iii) further by assuming instead
\begin{eqnarray*}
\lim_{k\to \infty}
\overline{\sE}^{(k)}_{j, \delta} (f,f) = {\sE}_{j, \delta} (f,f)
     \qquad \hbox{for every bounded measurable function } f \text{ on } B_j.
  \end{eqnarray*}
Then we can remove {\bf (A3)}${}^*$\,(ii). Note that
 $\sL_{j,\delta} f$ is bounded on $B_j$ for each $f\in \Lip_c(M)$ by \eqref{e:A1}.}
\end{remark}

\section{Weak convergence and discrete approximation}\label{Sect5}
\subsection{Weak convergence}
Let $X^{(k)}$ and $X$ be the symmetric Hunt processes associated with
$(\sE^{(k)}, \sF^{(k)})$ and $(\sE, \sF)$, respectively.

\begin{theorem}\label{fdc1}
Assume that
{\bf (A2)} holds and that $X$ is conservative.
We further assume that either
\eqref{e:A1a} of {\bf (A1)},
{\bf (A3)}--{\bf (A4)} hold, or {\bf (A3)}$^*$--\,{\bf (A4)}$^*$  hold.
Suppose $\varphi$ is in  $C^+_c(M)$. Then $\left\{X^{(k)}\right\}_{k \ge 1}$ with initial distribution
$\P^{(k)}_\varphi$ converge to $X$ with
initial distribution $\P_\varphi$ in the finite dimensional sense.
\end{theorem}

\proof
Without loss of generality, we assume  $\int_{M}\varphi(x)m(dx)=1$.
Let
$
P_t f (x) := \E_x[f(X_t)]$ and  $P^{(k)}_t g (x) := \E^{(k)}_x[g(X^{(k)}_t)]
$
be the contraction  semigroups on $L^2 (M, m)$ and $L^2 (V_k, m_k)$ respectively.
By Theorem \ref{Mosco1},  Theorem \ref{Mosco1--0}
and Theorem \ref{t:equi}, $E_kP^{(k)}_t \pi_k$ converges to $P_t$ strongly in $L^2 (M, m)$.
For any $l\geq 1$, $\{h_1, \cdots, h_l\}\subset L_b^2(M; m)$ and $0\leq t_1<t_2<\cdots < t_l$,
we have by Lemma \ref{l:sspi1} and the Markov property of $X^{(k)}$ and $X$ that
\begin{equation}\label{e:5.0} \lim_{k\to \infty}\E^{(k)}_{\varphi\cdot m_k} \left[ \pi_k h_1 (X^{(k)}_{t_1} ) \cdots \pi_k h_{l} (X^{(k)}_{t_{l}}) \right]
= \E_{\varphi \cdot m} \left[ h_1 (X_{t_1} ) \cdots h_{l} (X_{t_{l}}) \right].
\end{equation}
We fix $l \ge 1$.
Since $X$ is conservative, for any $\eps >0$, there is ball $B=B(x_0, r)$ so that
$\P_{\varphi \cdot m} (X_{t_j}\in B)>1-\eps$ for every $ j \in \{1, \dots l\}$.
By the strong $L^2$-convergence of $E_kP^{(k)}_{t_j} \pi_k \1_B$ to $P_{t_j}\1_B $
in $L^2 (M, m)$, we have
 \begin{equation}\label{e:5.2}  \lim_{k\to \infty} \P_{\varphi \cdot m_k}^{(k)}
\left( X^{(k)}_{t_j}\in B \right)>1-\eps  \quad \text{ for every } j \in \{1, \dots l\}.
\end{equation}
For any $\{f_1, \cdots, f_l\} \subset C_b (M)$,  since $E_k\pi_k f_j$ converges uniformly to $f_j$ on $\overline{B}$, from \eqref{e:5.0} we have
\begin{eqnarray}
 &&\lim_{k\to \infty}\E^{(k)}_{\varphi\cdot m_k}
 \left[f_1 (X^{(k)}_{t_1} ) \cdots   f_{l} (X^{(k)}_{t_{l}}):  \cap_{ j =1}^l\{X^{(k)}_{t_j}\in B\}\right]\nn\\
 &=& \lim_{k\to \infty}\E^{(k)}_{\varphi\cdot m_k}
 \left[\pi_k(f_1\1_B) (X^{(k)}_{t_1} ) \cdots  \pi_k (f_{l}\1_B) (X^{(k)}_{t_{l}})\right]\nn\\
 &=& \E_{\varphi\cdot m} \left[(f_1\1_B) (X_{t_1} ) \cdots (f_{j}\1_B) (X_{t_{j}}) \right]\nn\\
&=& \E_{\varphi\cdot m} \left[ f_1 (X_{t_1} ) \cdots f_{j} (X_{t_{j}}) : \cap_{ j =1}^l\{X_{t_j}\in B\}\right].\label{e:5.1}
\end{eqnarray}
We deduce the finite-dimensional convergence from \eqref{e:5.2} and \eqref{e:5.1}.\qed
\begin{definition}\label{def:s_separate} (\cite{EK})
Let $M$ be a metric space with metric $\rho$. A collection of function $S \subset C_b (M)$ is said to strongly separate points if for every $x \in M$ and $\delta > 0$, there exists a finite set $\{ h_1, \cdots , h_l \} \subset S$ such that
$$
\inf_{y: \rho(y,x) \ge \delta} \max_{1 \le i \le l} |h_i (y) - h_i (x) | > 0 .
$$
\end{definition}
\medskip
We can easily check that $\Lip^+_{c} (M)$  
strongly separates points in $M$.

\bigskip

\begin{theorem}\label{t:wc}
Assume that
{\bf (A1)}--{\bf (A2)} hold and that $X$ is conservative.
We further assume that either {\bf (A3)}--{\bf (A4)} hold, or {\bf (A3)}$^*$--\,{\bf (A4)}$^*$  hold.
Then, for
any $\varphi \in C^+_{c} (M)$,
 $\{(X^{(k)}, \, \P^{(k)}_\varphi); \, k
\ge 1\}$  converges weakly to $(X, \, \P_\varphi)$
 on $\bD_{M_\partial} [0, 1]$ equipped with the Skorohod topology.
\end{theorem}

\proof
First, note that,
by Proposition \ref{tight_aa},
for every $T>0$ and any $m \ge 1$ and
$\{g_1 , \cdots, g_m \} \subset \Lip^+_{c} (M)$, $\left\{( g_1 ,
\cdots, g_m)(X^{(k)})\right\}_{k \ge 1}$ restricted to $\{\zeta^{(k)}>T\}$
is tight in the  Skorohod space
$\bD_{\R^m} [0, T]$
with the initial distribution
$\P^{(k)}_{\varphi}$. Since $X$ is conservative, by \eqref{e:5.2}, for every $\eps >0$,
$$ \lim_{k\to \infty} \bP^{(k)}_{\varphi \cdot m_k} \left( \zeta^{(k)}>T \right) >1-\eps.
$$
So it follows from \cite[Theorem VI.3.21]{JS},
$\left\{( g_1 ,
\cdots, g_m)(X^{(k)})\right\}_{k \ge 1}$
is tight in the Skorohod space
$\bD_{\R^m} [0, T]$ with the initial distribution
$\P^{(k)}_{\varphi}$.
This together with Theorem \ref{fdc1} implies
the weak convergence of $\left\{( g_1 , \cdots,
g_m)(X^{(k)})\right\}_{k \ge 1}$ with initial distribution
$\P^{(k)}_\varphi$ to $( g_1 , \cdots, g_m)(X)$ with initial
distribution $\P_\varphi$. Since $\Lip^+_{c} (M)$  strongly
separates points in $M$, we have the desired result by Corollary
3.9.2 in \cite{EK}.    \qed

\subsection{Discrete approximation}\label{discrappr}
In this subsection, we give a general criteria for the approximation of
pure-jump process.

We give an extra condition on our approximating graphs.
\begin{description}
\item{\bf (AG.4)} There exists $n_0 \ge 1$ such that for every $j>n \ge n_0$ and
$x\in V_{2^j}$, there is some $y\in V_{2^n}$ so that $U_{2^j}(x)\subset U_{2^n}(y)$.
\end{description}
When $M=\bR^d$, the following approximation satisfies {\bf
(AG.1)}--{\bf (AG.4)}:
$V_k=k^{-1}\bZ^d$,
 $(x,y)\in B_k$ if and only if   $x,y\in k^{-1}\bZ^d$ with
  $\|x-y\|=k^{-1}$, and
for $V_k=\{x^{(k)}_{i}, i\geq 1\}$, $U_k(x^{(k)}_{i})=\prod_{k=1}^d [x^{(k)}_{i}-(2k)^{-1}, x^{(k)}_{i}+(2k)^{-1}]$.

Note that {\bf (AG.4)} is needed only this section.
Recall that $B_r = B(x_0, r)$ for $r>0$.

\begin{theorem}\label{t:latappr}
Let $j(x, y)$ be a non-negative measurable symmetric function on $M\times M$
 such that
$$ j(x, y) \leq M_0<\infty  \qquad \hbox{for every } x, y\in M \hbox{ with }
 \rho (x, y) \geq 1
$$
and for every compact set $K\subset M$,
$$ \lim_{j\to \infty} \sup_{x\in K} j(x, B_j^c) =0.
$$
Assume that the Dirichlet form $(\sE,\sF)$ determined by the jumping kernel
$j(x, dy):=j(x, y) m(dy)$ satisfies the conditions {\bf (A2)}--{\bf (A3)}.
Denote by $X$ the symmetric Hunt process associated with the
 regular Dirichlet form $(\sE,\sF)$ on $L^2(M,m)$, which we assume to be
 conservative.
  Let $(V_{2^k},B_{2^k}),\,k\in \bN$  be approximating graphs of $M$ and
$\{U_{2^k}(x)\}_{x\in V_{2^k}}$ be the associated partition satisfying
{\bf (AG.1)}--{\bf (AG.4)}.
Let
\begin{equation}\label{c^{(2^k)}aprx}
\sC^{(2^k)}(x,y):= \1_{\{\rho_{2^k}(x,y)\ge
4C_3/C_1\}} \, \frac {1}{m_{2^k}(x)m_{2^k}(y)}
\int_{U_{2^k}(x)} j(\xi , U_{2^k}(y)) m(d \xi)  \qquad x,y\in V_{2^k},
  \end{equation}
   where $m_{2^k}(x)=m(U_{2^k}(x))$ and $C_1,C_3$ are given
in (\ref{eq:equidis22}), (\ref{eq:equidis25}). Then
$(\sE^{(2^k)},\sF^{(2^k)})$ defined as in (\ref{e:Dfk1}) is a regular
Dirichlet form on  $L^2(V_{2^k},m_{2^k})$. Let $X^{(2^k)}$ be its associated
continuous time Markov chain on $V_{2^k}$. Then, for any positive
function
$\varphi \in C^+_{c} (M)$,
$\{(X^{(2^k)}, \,
\P^{(2^k)}_\varphi); \, k \ge 1\}$  converges weakly to $(X, \,
\P_\varphi)$
on $\bD_{M_\partial} [0, 1]$ equipped with the Skorohod topology.
\end{theorem}

\proof  For notational simplicity, in this proof we write $k$ for $2^k$.
In view of  Theorem \ref{t:wc}, it is enough to show {\bf (A1)} and {\bf (A4)} hold.
 For $\rho_k(x,y)\ge 4C_3/C_1$ and $\xi\in U_k(x), \eta\in
U_k(y)$,  we  have by \eqref{eq:equidis22}--\eqref{eq:equidis25}
and the triangle
inequality that
$\rho(x,y)\ge C_1\rho_k(x,y)/k\ge 4C_3/k$,
\begin{equation}\label{e:5.4}
 |\rho(\xi,\eta)- \rho (x, y)| \leq \rho (x,  \xi)+\rho (\eta, y)
\leq  C_3/k + C_3/k= 2C_3/k
 \end{equation}
and so
\begin{equation}\label{e:5.5}
\frac{C_1}{2} \frac{\rho_k(x, y)}{k} \leq
 \rho (x, y)/2 \leq \rho (\xi , \eta) \leq 3\rho (x, y)/2
 \leq \frac{3C_2}{2} \frac{\rho_k(x, y)}{k} .
\end{equation}
Take a compact set $K \subset M$ and $K_1:=\{x \in M: \rho(x,K) \le 1  \}$.
Then by \eqref{e:5.5}
\begin{eqnarray*}
&&\sup_{k\in \bN}\sup_{x \in K\cap V_k} \sum_{y\in V_k} \sC^{(k)}(x,y)
\Big( \frac {\rho_k(x,y)}k  \wedge 1\Big)^2 m_{k}(y) \\
&=&\sup_{k\in \bN}\sup_{x \in K\cap V_k} \sum_{y\in V_k}
\Big(\frac {\rho_k(x,y)}k   \wedge 1\Big)^2 \,  \1_{\{\rho_k(x,y)\ge 4C_3/C_1\}}
 \frac {1}{m_k(x)} \int_{U_k(x)} j(\xi, U_k(y))
 m(d\xi) \\
 &=&\sup_{k\in \bN}\sup_{x \in K\cap V_k} \sum_{y\in V_k}
 \frac {1}{m_k(x)} \int_{U_k(x)}\int_{U_k(y)} \Big(\frac {\rho_k(x,y)}k   \wedge 1\Big)^2  \1_{\{\rho_k(x,y)\ge 4C_3/C_1\}} \,  j(\xi, d\eta)
 m(d\xi) \\
 &\le & c \sup_{k\in \bN}\sup_{x \in K\cap V_k} \sum_{y\in V_k}
 \frac {1}{m_k(x)} \int_{U_k(x)} \left(\sup_{\xi \in U_k(x)}  \int_{U_k(y)}
(\rho(\xi,\eta)^2 \wedge 1)
  j(\xi, d\eta) \right)
 m(d\xi) \\
&\le &c\sup_{k\in \bN}\sup_{\xi \in K_1} \sum_{y\in V_k}
\int_{U_k(y)}(\rho(\xi,\eta)^2 \wedge 1) j(\xi,d\eta) \\
&\le&  c\sup_{\xi \in K_1}\int_{M}(\rho(\xi,\eta)^2 \wedge 1) j(\xi,d\eta)\le C_K
\end{eqnarray*}
by {\bf (A2)}\,(iii). This
proves \eqref{e:A1a}  of {\bf (A1)}.

 By \eqref{e:5.4}, for $k\geq 2C_3$ and $x, y\in V_k$ with $\rho_k (x, y) \geq 2$,
 $$ \rho (\xi, \eta )  \geq \rho (x, y)- 2C_3/k\geq 1 \qquad \hbox{for } \xi \in U_k(x) \hbox{ and }
  \eta \in U_k(y).
  $$
So for each $k\geq 2C_3$, $j\geq 1$ and $x\in \overline B_j\cap V_k$, $y\in
B_{j+2}^c\cap   V_k$,
$$\sC^{(k)}(x, y) \leq  \frac {1}{m_k(x)m_k(y)}
\int_{U_k(x)\times U_k(y)} j(\xi , \eta) m(d \xi) m(d\eta) \leq M,
$$
which establishes \eqref{e:A1b} of {\bf (A1)}.

By definition of $C^{(k)}(\cdot,\cdot)$, (\ref{e:A4.1}) clearly holds.
For any compact set $K\subset M$ with $K_1:=\{x\in M: \rho (x, K)\leq 1\}$, we have
\begin{eqnarray*}
\lim_{j\to \infty} \sup_{k\geq 1} \sup_{x\in K} \int_{B_j^c}  \sC^{(k)} (x, y) m(dy)
&\leq & \lim_{j\to \infty} \sup_{x\in K_1}
\int_{ B_j^c}  j (\xi , y)   m(dy)=0,
\end{eqnarray*}
so (\ref{e:A4.2}) holds.

On the other hand by {\bf (A1)}, for any
$f \in L^2_b(M)$
with $\|f\|_\infty \le M_1$, $j \ge 1$ and $\delta>0$,
\begin{eqnarray}
 && \Big| \overline\sE^{(k)}_{j,\delta} (f,f)^{1/2} -
    \overline\sE^{(k)}_{j,\delta} (E_k \pi_kf,  E_k \pi_k f)^{1/2}\Big| \nonumber \\
   & \leq & \overline\sE^{(k)}_{j,\delta} (E_k \pi_kf-f,  E_k \pi_k f-f)^{1/2}\nonumber  \\
    &\leq & \left(
   2 \int_{B_j} (f(x)-E_k \pi_k f(x))^2 \left(\int_{B_j}  \sC^{(k)} (x, y)  \1_{\{\rho (x, y) > \delta\}} m(dy)\right) m(dy)\right)^{1/2} \nonumber  \\
   &\leq & c(j, \delta) \, \|  f-E_k \pi_kf\|_2 , \label{e:5.9}
\end{eqnarray}
which goes to 0 as $k\to \infty$ by Lemma
\ref{l:sspi1}(iv).
Note that
for large $k$ and small $\delta$,
       \begin{eqnarray}
 &&\overline\sE^{(k)}_{j,\delta} (E_k \pi_kf,  E_k \pi_k f)=  \int_{B_j\times B_j} (E_k\pi_k f(x)-E_k \pi_k f(y))^2 \sC^{(k)} (x, y) \1_{\{\rho (x, y) > \delta\}} m(dx)m(dy)   \nn\\
    &=& \frac12
 \sum_{(z,w) \in V_k\times V_k} (\pi_k f(z)-\pi_k f(w))^2   \, \frac {1}{m_{k}(z)m_{k}(w)}
\int_{U_{k}(z)} j(\xi , U_{k}(w)) m(d \xi) \times \nn \\
&&\quad \times \int_{(B_j\times B_j) \cap (U_k(z) \times U_k(w))}
  \1_{\{\rho (x, y) > \delta\}} m(dx)m(dy)   \label{e:newq1}
     \end{eqnarray}
  and
  \begin{eqnarray}
  &&  \sE_{j, \delta} (E_k \pi_kf,  E_k \pi_k f)
 = \frac12
  \int_{B_j\times B_j} (E_k\pi_k f(x)-E_k \pi_k f(y))^2 j (x, y) \1_{\{\rho (x, y) > \delta\}} m(dx)m(dy)\nn\\
    &=& \frac12\sum_{(z,w) \in V_k\times V_k}  (\pi_k f(z)- \pi_k f(w))^2\int_{(B_j\times B_j) \cap (U_k(z) \times U_k(w))}   j (x, y) \1_{\{\rho (x, y) > \delta\}} m(dx)m(dy)\label{e:newq2}
       \end{eqnarray}
Since summands in \eqref{e:newq1} and \eqref{e:newq2} are same except the case $\rho(x,y)$ small and $y$ is near the boundary of $B_j$, it is easy to see that
      there exists $k_0 =k_0(\delta)>0$ and $c>0$ such that  for $k \ge k_0$,
\begin{eqnarray*}
 &&   \Big|
    \overline\sE^{(k)}_{j,\delta} (E_k \pi_kf,  E_k \pi_k f)-\sE_{j, \delta} (E_k \pi_kf,  E_k \pi_k f)\Big|  \\
 &\le &
  2\int_{B_j\times \{j-c\frac1{k}<\rho(y, x_0)<j+ c\frac1{k}\}} (E_k\pi_k f(x)-E_k \pi_k f(y))^2 j(x,y) \1_{\{\rho (x, y) > \delta -c\frac1{k}\}} m(dx)m(dy)\\
&&  +\int_{B_{j+1}\times B_{j+1}} (E_k\pi_k f(x)-E_k \pi_k f(y))^2 j (x, y) \1_{\{\delta +c\frac1{k}>\rho (x, y) > \delta-c\frac1{k}\}} m(dx)m(dy)\\
  &\le&2(2M_1)^2\int_{B_j\times \{j-c\frac1{k}<\rho(y, x_0)<j+ c\frac1{k}\}}  j(x,y) \1_{\{\rho (x, y) > \delta -c\frac1{k}\}} m(dx)m(dy)\\
&&  +(2M_1)^2\int_{B_{j+1}\times B_{j+1}} j (x, y) \1_{\{\delta +c\frac1{k}>\rho (x, y) > \delta-c\frac1{k}\}} m(dx)m(dy),
     \end{eqnarray*}
    which goes to zero as $k$ goes to $\infty$.
    Therefore
    \begin{eqnarray}
 && \lim_{k\to\infty}  \Big|
    \overline\sE^{(k)}_{j,\delta} (E_k \pi_kf,  E_k \pi_k f)-\sE_{j, \delta} (f,f)\Big|  \nonumber\\
    &\le &
     c(M_1,j,\delta, f)
    \lim_{k\to\infty}  \Big|
    \overline\sE^{(k)}_{j,\delta} (E_k \pi_kf,  E_k \pi_k f)^{1/2}-\sE_{j, \delta} (f,f)^{1/2}\Big|  \nonumber\\
 &\le&
     c(M_1,j,\delta, f)\lim_{k\to\infty}  \Big|
    \sE_{j, \delta} (E_k \pi_kf,  E_k \pi_k f)^{1/2}-\sE_{j, \delta} (f,f)^{1/2}\Big| \nonumber \\
    &\leq &
     c(M_1,j,\delta, f)\lim_{k\to\infty} \sE_{j, \delta} (E_k \pi_kf-f,  E_k \pi_k f-f)^{1/2}\nonumber  \\
    &\leq &
     c(M_1,j,\delta, f)\lim_{k\to \infty}   \left( \int_{B_j\times B_j} ((f-E_k\pi_k f)(x)-(f-E_k \pi_k )f(y))^2 j (x, y) \1_{\{\rho (x, y) > \delta\}} m(dx)m(dy)\right)^{1/2} \nonumber\\
    &\leq &
     c(M_1,j,\delta, f) \lim_{k\to \infty} \left( \int_{B_j} (f(x)-E_k \pi_k f(x))^2
    \left( \int_{B_j} j (x, y) \1_{\{\rho (x, y) > \delta\}} m(dy\right) m(dx)\right)^{1/2}\nonumber \\
    &\le &
    c(M_1, j, \delta, f)
     \, \lim_{k\to \infty} \|  f-E_k \pi_kf\|_2  =0.\label{eq:pandeij}
    \end{eqnarray}
This combined with
\eqref{e:5.9} shows that
$
 \lim_{k\to \infty}\overline\sE^{(k)}_{j,\delta} (f,f)=\sE_{j, \delta} (f,f)$
 for any $f \in L^2_b(M)$.

 The monotonicity property of
 {\bf (A4)}(ii) (with $2^k$ instead of $k$)
  is an immediate consequence of {\bf (AG.4)}
 and \eqref{c^{(2^k)}aprx}.
So we have established {\bf (A4)}. \qed

\begin{remark} {\rm
For any $f \in L^2_b(M)$ with
$\|f\|_\infty \le M_1$, $j \ge 1$ and $\delta>0$,
computing similarly to \eqref{eq:pandeij}, we have
\begin{eqnarray}
 \Big| \| \overline \sL^{(k)}_{j,\delta} f\|_{2,B_j}-
 \|\overline \sL^{(k)}_{j,\delta} E_k \pi_kf\|_{2,B_j}  \Big|
    \leq  \Big| \| \overline \sL^{(k)}_{j,\delta} (f- E_k \pi_kf)\|_{2,B_j}  \Big|    \leq  c(j, \delta) \, \|  f-E_k \pi_kf\|_2.  \nonumber
\end{eqnarray}
which goes to 0 as $k\to \infty$ by Lemma \eqref{l:sspi1} (iv).
Moreover, by Lemma \ref{l:bas_lf} (i) (ii),
   \begin{eqnarray*}
  \lim_{k\to\infty} \Big| \|  \sL_{j,\delta} f\|_{2,B_j}-
 \|\sL_{j,\delta} E_k \pi_kf\|_{2,B_j} \Big|
   & \leq &\Big| \| \sL_{j,\delta} (f- E_k \pi_kf)\|_{2,B_j}  \Big| \nonumber  \\
 &\le & \lim_{k\to\infty}  c(j, \delta) \, \|  f-E_k \pi_kf\|_2 =0
\end{eqnarray*}

Thus, to show {\bf (A4)${}^*$}\,(ii),  it is enough to show that
\begin{equation}\label{e:dpp}\limsup_{k\to\infty} |\|\overline \sL^{(k)}_{j,\delta} E_k \pi_kf\|^2_{2,B_j}  -\|\sL_{j,\delta} E_k \pi_kf\|^2_{2,B_j}|=0.
\end{equation}
}
\end{remark}

\section{Tightness and weak convergence under convergence-in-measure topology}\label{Sect6}

In some of the applications, we need the convergence-in-measure topology
on  $\bD_{M_\partial} [0, 1]$ and on $\bD_{M} [0, 1]$ , which was introduced in \cite{DM}
and is weaker than the Skorohod topology.
This convergence-in-measure topology is also called pseudo-path
topology in literature, see \cite[Lemma 1]{MZ}.

Let $\lambda$ be the Lebesgue measure on $[0,1]$. For a $M_\partial$-valued Borel function on $[0,1]$,
the pseudo-path of $w$ is a probability law on $[0,1]\times M_\partial$: the image measure of
$\lambda$ under the mapping $t\mapsto (t, w(t))$. Denote by $\Psi$ the mapping which associates to
a path $w$ its pseudo-path, which identifies two paths if and only if they are equal $\lambda$-a.e.
on $[0,1]$. In particular, $\Psi$ is one-to-one on $D_{M_\partial}[0,1]$
and embeds it into the compact space of all probability measures on the compact space
$[0,1]\times M_\partial$. Meyer
gave the name of the pseudo-path topology to the induced topology
on $D_{M_\partial}[0,1]$.
(See \cite[chapter IV, n 40-46]{DM} for more details.) 
Theorem 5 of \cite{MZ} tells us that if the law of $\{X^{k}, k\geq 1\}$ is tight in $D_{M_\partial}[0, 1]$ equipped with pseudo-path topology, then there is a subsequence $\{n_k\}$ and a subset $A$ of $[0, 1]$ having zero Lebesgue
measure so that $X^{n_k}$ convergence in finite dimensional distribution  on $[0, 1]\setminus A$.

Tightness of stochastic processes on $\bD_{M_\partial} [0, 1]$
(respectively, on $\bD_{M} [0, 1]$) equipped with
the convergence-in-measure topology is closely related to the
number of crossing between two disjoint sets by the stochastic
processes (see \cite{MZ}). The latter has been investigated in \cite{CFS,LZ}.

\begin{proposition}\label{T:5.1}
Assume that {\bf (A.2)}, {\bf (A.3)} and
{\bf (A.4)}{\rm (i)(iii)} hold. Then  for every $\varphi \in C^+_c(M)$, the law
$\{\P^{(k)}_\varphi, k\geq 1\}$ is tight on
$\bD_{M_\partial} [0, 1]$ equipped with the convergence-in-measure
topology.
\end{proposition}
\proof Let $D_1$ and $D_2$ be two relatively compact open subsets
in $M$ with disjoint closure.
By  {\bf (A.3)}, there is some $f\in \Lip_c(M)\subset \sF$
so that $f=1$ in an open
 neighborhood of $\overline D_2$ and $f=0$ in an open neighborhood
 of $\overline D_1$. Then for $k$ sufficiently large,
 $\pi_k f=1$ on $V_k\cap D_2$ and $\pi_k f=0$ on $V_k\cap D_1$.
Let $N^{(k)}$ be the number of crossings by $X^{(k)}$ from $D_1$
into $D_2$. By \cite[Theorem in page 69]{CFS},
if $g\in \sF^{(k)}$ such that $g=1$ on $D_2\cap V_k$ and $g=0$ on $D_1\cap V_k$,
then
\begin{equation}
 \E^{(k)}_{\varphi \cdot m_k} [ N^{(k)}]
\leq 2 \| \varphi \|_\infty\,   \sE^{(k)} (g, g).
\end{equation}
It follows from Lemma \ref{l:fcon1} that
$$\sup_{k\geq 1}  \E^{(k)}_{\varphi \cdot m_k} [ N^{(k)}] <\infty.
$$
Since the above holds for every pair of relatively compact open subsets
in $M$ with disjoint closure, we conclude by \cite[Theorem 2]{MZ}
and a diagonal selection procedure that the law
$\{\P^{(k)}_\varphi, k\geq 1\}$ is tight on
$\bD_{M_\partial} [0, 1]$ equipped with the convergence-in-measure
topology.  \qed

\begin{theorem}\label{T:5.2}
Assume that either  \eqref{e:A1a} of {\bf (A1)} and {\bf (A2)}--{\bf (A4)} hold,  or
 {\bf (A.2)}, {\bf (A.3)${}^*$}
and
{\bf (A.4)${}^*$} hold. Then  for every $\varphi \in C^+_c(M)$,
 $\{(X^{(k)},$ $ \P^{(k)}_\varphi); \, k
\ge 1\}$  converges weakly to $(X, \, \P_\varphi)$
 on $\bD_{M_\partial} [0, 1]$ equipped with the convergence-in-measure topology,
 where $X$ is the Hunt process  associated with $(\sE, \sF)$.
\end{theorem}

\proof First, note that conditions {\bf (A.3)${}^*$} and
{\bf (A.4)${}^*$} are stronger than conditions {\bf (A.3)} and
{\bf (A.4)}(i)(iii). So, by Proposition \ref{T:5.1}, for any subsequence $\{n_k; k\ge 1\}$, there exists a
sub-subsequence $\{n_k'; k\ge 1\}$ such that  $\{(X^{(n_k')}, \, \P^{(n_k')}_\varphi); k\ge 1\}$
converges weakly on $\bD_{M_\partial} [0, 1]$ equipped with the convergence-in-measure topology
to a law of say $\tilde {\mathbb P}$. 
Thus by \cite[Theorem 5]{MZ}, we may assume without
loss of generality that there is a subset $A\subset [0, 1]$
of zero Lebesgue measure so that 
$\{(X^{(n_k')}, \, \P^{(n_k')}_\varphi); k\ge 1\}$ converges
in finite dimension over the time interval $[0, 1]\setminus A$ to that 
of $\tilde {\mathbb P}$.
Let
$
P_t f (x) := \E_x[f(X_t)]$ and  $P^{(k)}_t g (x) := \E^{(k)}_x[g(X^{(k)}_t)].
$
By Theorem \ref{Mosco1} or Theorem \ref{Mosco1--0},
we know that
$(\sE^{(k)}, \sF^{(k)})$ is Mosco convergent to $(\sE, \sF)$.
 So
 by Theorem \ref{t:equi} (ii), $E_kP^{(k)}_t \pi_k f$ converges to $P_tf$  in $L^2 (M, m)$. 
This implies by the Markov property that,  for any $l\geq 1$, $\{h_1, \cdots, h_l\}\subset C_c(M; m)$ and $0\leq t_1<t_2<\cdots < t_l$,
$$ \lim_{k\to \infty}\E^{(k)}_{\varphi\cdot m_k} \left[ \pi_k h_1 (X^{(k)}_{t_1} ) \cdots \pi_k h_{l} (X^{(k)}_{t_{l}}) \right]
= \E_{\varphi \cdot m} \left[ h_1 (X_{t_1} ) \cdots h_{l} (X_{t_{l}}) \right].
$$
Thus the finite dimensional distribution under $\tilde {\mathbb P}$ over
the time interval $[0, 1]\setminus A$ is the same as that of 
$(X, \, \P_\varphi)$. Since both laws $\tilde {\mathbb P}$ and $\bP_\varphi$
are carried on $\bD_{M_\partial}[0, 1]$, 
 it follows 
that $\tilde {\mathbb P}$ has the same distribution as the law of  $(X, \, \P_\varphi)$.
Since this holds for any subsequence $\{n_k; k\ge 1\}$, we obtain the desired result. 
\qed

\section{Application to random walk in random conductance}\label{applisec}

In this section, we 
present  application of Theorem \ref{Mosco1--0} to the scaling limit of
some random walk in random conductance.

Throughout this subsection, $M=\mathbb R^d$
and $m$ be a $d$-dimensional Lebesgue measure. Also, let $V_k=k^{-1}\mathbb Z^d$ and
$m_k(x)=k^{-d}$ for all $x\in V_k$.

Let
$j(x,y)$ be  a symmetric non-negative continuous function
of $x$ and $y$ on $\bR^d \times \bR^d \setminus \wh d$ such that
 there
exist $\alpha, \beta \in (0,2)$, $\alpha > \beta $ and positive
$\kappa_1, \kappa_2$ such that
 \begin{equation}\label{e:3.5}
\kappa_1 |y-x|^{-d-\beta} \leq j(x,y) \leq \kappa_2
|y-x|^{-d-\alpha} \quad \hbox{for } |y-x|<1
\end{equation}
and
 \begin{equation}\label{e:3.6}
\sup_{(x,y) \in \bR^d \times \bR^d \atop |y-x|\ge 1}j(x,y)  \le \kappa_0 < \infty    \quad \text{and }\quad \sup_{x \in \R^d} \int_{\{|y-x|\ge 1\}} j(x,y) m(dx) < \infty.
\end{equation}
Set the Dirichlet form $(\sE,\sF)$ which is defined by (\ref{e:Df}) with $J(dx,dy)=j(x,y)m(dx)m(dy)$ where $j(x,y)$ defined in \eqref{e:3.5}--\eqref{e:3.6}. Finally we assume {\bf (A3)} is true. i.e., $\Lip_c(M)$ is dense in $(\sF,
\sE(\cdot,\cdot)+\|\cdot\|_2^2)$.
Then,  by \cite[Propostion 2.2]{CK2} and its proof,
the Dirichlet form $(\sE,
\sF)$ is regular on $\R^d$ and so it associates a Hunt process $X$
starting from quasi-everywhere in $\R^d$. Moreover $X$ is conservative since \eqref{e:cons2} holds.

\begin{proposition}\label{arrc}
{\rm (i)}
Suppose $d \ge 2$.
Let $\{\xi_{x,y}\}_{x,y\in \bZ^d, x\ne y}$ be
i.i.d. on $(\Omega,
\sA, {\bf P})$ such that
$0\le\xi_{x,y}$,
${\bf E}[\xi_{x,y}] =1$ and $\mbox{Var }(\xi_{x,y})<\infty$.
Let
\begin{equation}\label{eq:eneRCM}
\sC^{(k)}(x,y):= \xi_{kx, ky} j(x,y)
~~~\mbox{for }~~x,y\in V_k.\end{equation}
Let $(\sE^{(k)},\sF^{(k)})$ be defined as in (\ref{e:Dfk1}) and
define the Markov chain corresponding to
$\sE^{(k)}$ by $X^{(k)}_t$.
Let $X$ be the Hunt process corresponding to $(\sE,\sF)$
which is defined by (\ref{e:Df})
 with $J(dx,dy)=j(x,y)m(dx)m(dy)$ where $j(x,y)$ defined in \eqref{e:3.5}--\eqref{e:3.6}.
Define
$T^{(k)}_t$ and $T_t$ as the semigroups corresponding to $X^{(k)}$ and $X$ respectively.
Then,
$E_k T^{(k)}_t \pi_k \to T_t$ strongly in $L^2(\bR^d,m)$ ${\bf P}$-a.s. and the convergence
is uniform in any finite interval of $t \ge 0$.
Moreover, $(X^{(k)}, \P^{(k)}_\varphi)$ converges weakly to $(X, \P_\varphi)$ on $\bD_{M_\partial} [0, 1]$ equipped with convergence-in-measure topology  ${\bf P}$-a.s..

{\rm (ii)} Assume  further that $0\le\xi_{x, y}\le C$
${\bf P}$-a.s.
for some deterministic constant $C>0$. Then for any positive function
$\varphi \in C_{c} (M)$, $\{(X^{(k)}, \, \P^{(k)}_\varphi); \, k
\ge 1\}$  converges weakly to $(X, \, \P_\varphi)$
on $\bD_{M_\partial} [0, 1]$ equipped with the Skorohod topology  ${\bf P}$-a.s..
\end{proposition}

\proof (i) Note first that since, by \eqref{e:3.6}
$$
{\bf E}[ \sum_{y\in V_k} \sC^{(k)}(x,y) m_{k}(y)] \le \kappa_2 \sum_{y\in V_k, |x-y| <1}
k^{-d}|x-y|^{-d-\alpha}+\sum_{y\in V_k, |x-y| \ge 1}
k^{-d}j(x,y)<\infty,$$
we have $ \sum_{y\in V_k} \sC^{(k)}(x,y) m_{k}(y)
<\infty$ ${\bf P}$-a.s., so (\ref{e:2.10}) holds. Thus, by Theorem
\ref{T:2.4}, $(\sE^{(k)},\sF^{(k)})$ is a regular Dirichlet form.
In order to prove the first assertion of (i),
by Theorem \ref{Mosco1--0},
 Theorem \ref{fdc1}
 and Theorem \ref{t:equi}, it is enough to prove {\bf (A2)}, {\bf (A3)${}^*$} and
{\bf (A4)${}^*$} ${\bf P}$-a.s..
Recall that we assume {\bf (A3)}. Moreover,  {\bf (A3)${}^*$}(ii) is true by the continuity of $j(x,y)$. Furthermore, by symmetry of
 $j(x,y)$ and \eqref{e:3.5}--\eqref{e:3.6}, one can easily see that
{\bf (A2)} is true. So, we
will prove {\bf (A4)${}^*$} below.

We first show (\ref{e:A4.1}). Let $\eta \le 1$. Note that, by \eqref{e:3.5}
\[\int\int_{\{(x, y)\in K\times K: |x- y|\leq  \eta\}} |x- y|^2 \sC^{(k)}(x,y)
   m(dx)m(dy)
 \le \kappa_2 k^{-2d} \sum_{{x,y\in V_k\cap K}\atop{|x-y|\le \eta}}\frac{|x-y|^2\xi_{kx,ky}}{|x-y|^{d+\alpha}}
   =:\kappa_2 k^{-2d}H_k.\]
Since $|x-y|\ge k^{-1}$ when $x\ne y$, setting $2-\alpha=\eps$,
\[\mbox{Var }(H_k)=
\sum_{{x,y\in V_k\cap K}\atop{|x-y|\le \eta}}|x-y|^{2(2-d-\alpha)}
\mbox{Var }(\xi_{kx,ky})
\le c_1k^{3d}\cdot k^{-2d}\sum_{{x,y\in V_k\cap K}\atop{|x-y|\le \eta}}|x-y|^{-d+2 \eps}
\le c_2k^{3d}m(K)\eta^{\eps}.\]
So,
\[{\bf P}(k^{-2d}\left|H_k-{\bf E}[H_k]\right|\ge \eta^{\eps/2})
\le \kappa^2_2
\frac{\mbox{Var }(H_k)}{k^{4d}\eta^{\eps}}
\le \frac{c_3}{k^{d}},\]
and using the Borel-Cantelli Lemma, we have $\limsup_kk^{-2d}|H_k-{\bf E}[H_k]|\le \eta^{\eps/2}$
${\bf P}$-a.s., so
\[\lim_{\eta\to 0}\limsup_{k\to\infty}k^{-2d}|H_k-{\bf E}[H_k]|=0.\]
On the other hand, by \eqref{e:3.5}
\begin{eqnarray*}
\limsup_{k\to\infty}k^{-2d}{\bf E}[H_k]
&\le & \kappa_2
\limsup_{k\to\infty}k^{-2d}\sum_{{x,y\in V_k\cap K}\atop{|x-y|\le \eta}}|x-y|^{(2-d-\alpha)}
{\bf E}[\xi_{kx,ky}]\\
&=&
\kappa_2\limsup_{k\to\infty}k^{-2d}\sum_{{x,y\in V_k\cap K}\atop{|x-y|\le \eta}}|x-y|^{2-d-\alpha}
\,\le\, cm(K)\eta^{(2-\alpha)/2},
\end{eqnarray*}
 which
vanishes when $\eta\to 0$, so we obtain (\ref{e:A4.1}) ${\bf P}$-a.s..

We next show (\ref{e:A4.2}).
Note that
\[\int_K \int_{B_j^c} \sC^{(k)} (x, y) m(dx)m(dy)
=k^{-2d}\sum_{{y\in V_k\cap K}}
   \sum_{x\in V_k\cap B_j^c}{\xi_{kx,ky}}j(x,y)
   =:k^{-2d}H'_k.\]
Then,
for $j \ge j_0$ where $K \subset B_{j_0-1}$, by \eqref{e:3.6} we have
\begin{eqnarray*}
k^{-2d}\mbox{Var }(H'_k)&=&k^{-2d}
\sum_{{x\in V_k\cap B_j^c}\atop{y\in V_k\cap K}}
\mbox{Var }(\xi_{kx,ky}) j(x,y)^2\\
&\le& c k^{-2d}
\sum_{{x\in V_k\cap B_j^c}\atop{y\in V_k\cap K}}
 j(x,y)
\le \,c \,k^{-2d}
\sum_{{x\in V_k : |x-y| > j-j_0}\atop{y\in V_k\cap K}}
 j(x,y)\,=: \,c \,a_j^k
\end{eqnarray*}
Thus,
\[{\bf P}(k^{-2d}(a^k_j)^{-1/2}\left|H'_k-{\bf E}[H'_k]\right|\ge 1)
\le
\frac{\mbox{Var }(H'_k)}{k^{4d}a^k_j}
\le \frac{c}{k^{2d}}      ,\]
and using the Borel-Cantelli Lemma, we have $\limsup_k k^{-2d}(a^k_j)^{-1/2}|H'_k-{\bf E}[H'_k]|\le 1$
${\bf P}$-a.s..
Since $a^k_j$ converges to
$$ a_j:=\int_K \int_{\{|x-y|>j-j_0\}} j(x,y)m(dx)m(dy) \in (0,  \infty)$$
 by continuity of $j(x,y)$ and \eqref{e:3.6}, we have
\[\lim_{j\to \infty}\limsup_{k\to\infty}k^{-2d}|H'_k-{\bf E}[H'_k]|=
(\limsup_{k\to\infty} k^{-2d}(a^k_j)^{-1/2}|H'_k-{\bf E}[H'_k]|) \lim_{j\to \infty}\sqrt{a_j}
\le \lim_{j\to \infty}\sqrt{a_j}=0
.\]
In the last equality above, we have used \eqref{e:3.6}.
On the other hand, by similar computation we have
\begin{eqnarray*}
\lim_{j\to \infty}\limsup_{k\to\infty}k^{-2d}{\bf E}[H'_k]\le c\lim_{j\to \infty} a_j=0\end{eqnarray*}
We have proved (\ref{e:A4.2}).

For the remainder part of the proof, we fix $\delta,  j>0$.
We now show {\bf (A4)${}^*$}\,(iii).

Let $h$ be a bounded  and continuous function in $B_j\times B_j$.
By the continuity and boundedness of $h(x,y)$ and $j(x,y)$ on $B_j \times B_j \setminus \wh d$,
 we have
\begin{equation}\label{eq:nirbw}   \lim_{k\to \infty}
k^{-2d}
\sum_{{x, y\in V_k\cap B_j}\atop{|x-y|>\delta}}
 {h(x,y)}j(x,y)
= \int_{   B_j\times B_j} h(x,y) \, \1_{\{ |x-y|>\delta\}}\,  j(x,y)m(dx)m(dy),\end{equation}
so it is enough to show
\begin{equation}\label{q:dhealed}
\lim_{k\to \infty}
k^{-2d}
\sum_{{x, y\in V_k\cap B_j}\atop{|x-y|>\delta}}
 {h(x,y)(\xi_{kx,ky}-1)}j(x,y)=0~~~~\mbox{{\bf P}-a.s.}. \end{equation}
Using \eqref{e:3.5}--\eqref{e:3.6}, we have,
\begin{eqnarray*}
&&{\bf P}\Big(k^{-2d}\,\Big|
\sum_{{x, y\in V_k\cap B_j}\atop{|x-y|>\delta}}
  {h(x,y)(\xi_{kx,ky}-1)}j(x,y)\Big|>\eps^{1/2}\Big)\\
 &\le &
 c_1   \frac{1}{k^{4d}\eps} \mbox{Var }\Big(\sum_{{x, y\in V_k\cap B_j}\atop{|x-y|>\delta}}
  {h(x,y)(\xi_{kx,ky}-1)}j(x,y)\Big)\\
  &\le &
 c_2   \frac{1}{k^{2d}\eps}\mbox{Var } \left(\xi_{kx,ky}\right) \Big(\frac{1}{k^{2d}} \sum_{{x,y\in V_k\cap B_j}\atop{|x-y|>\delta}}h(x,y)^2 |x-y|^{-2d-2\alpha} \Big)
 \le  \frac{c_{\delta,j}}{k^{2d}\eps},
 \end{eqnarray*}
so using the Borel-Cantelli Lemma, computing similarly as before, we obtain (\ref{q:dhealed}).

Lastly, we show {\bf (A4)${}^*$}\,(ii). Fix $f \in \Lip_c(M)$. Note that
\begin{eqnarray*}
\overline \sL^{(k)}_{j,\delta}f(x)&=&\frac 1{k^d}\sum_{{y\in V_k\cap B_j}\atop{|x-y|>\delta}}
{(f(y)-f(x))}j(x,y)+
\frac 1{k^d}\sum_{{y\in V_k\cap B_j}\atop{|x-y|>\delta}}
{(\xi_{kx,ky}-1)(f(y)-f(x))}j(x,y)\\
&=:&I^{(k)}_1(x)+I^{(k)}_2(x).
\end{eqnarray*}

One can easily see that $\|I^{(k)}_1-\sL_{j,\delta} f\|_2\to 0$ as $k\to \infty$.
Indeed, by the continuity and boundedness of $j$ and $f$, it is clear that $\lim_{k \to \infty}I^{(k)}_1(x)=\sL_{j,\delta} f(x)$ for all $x$ and
$|I^{(k)}_1(x)|\le C$ for large $C$.
Thus the bounded convergence theorem can be applied.
So all we need is to show $\|I^{(k)}_2\|_2\to 0$\, {\bf P}-a.s. as $k\to \infty$.
Since
\begin{eqnarray*}
{\bf E}[\|I^{(k)}_2\|_2^2]
&=& k^{-2d}{\bf E}\Big[  \int_{B_j}      \Big(  \sum_{{y\in V_k\cap B_j}\atop{|x-y|>\delta}}
{(\xi_{kx,ky}-1)(f(y)-f(x))}j(x,y)      \Big)^2 m(dx) \Big]\\
&=& k^{-2d}  \int_{B_j}    \sum_{{y\in V_k\cap B_j}\atop{|x-y|>\delta}}  {(f(x)-f(y))^2
\mbox{Var }(\xi_{kx,ky})} j(x,y)^2m(dx) \\
&=& c k^{-d}  \int_{B_j}    \sum_{{y\in V_k\cap B_j}\atop{|x-y|>\delta}}  {(f(x)-f(y))^2}
 j(x,y)^2 m_k(y)  m(dx)
\,\le\, c_{f,\delta,j}k^{-d},
\end{eqnarray*}
 computing similarly as before,
\bee\label{eq:consrst}
{\bf P}(
\|I^{(k)}_2\|_2^2>\eps)\le \eps^{-1}{\bf E}[\|I^{(k)}_2\|_2^2]\le
\frac{c_{f,\delta,j}}{\eps k^d}.\eee
So using the Borel-Cantelli Lemma,  $\|I^{(k)}_2\|_2\to 0$\, {\bf P}-a.s. for $d\ge 2$.
The weak convergence follows from
Theorem \ref{T:5.2}.

(ii)
Using \eqref{e:3.5}--\eqref{e:3.6}, it is easy to show that {\bf (A1)} holds ${\bf P}$-a.s., and
 $X$ is conservative.
 Thus,
  by Theorem \ref{Mosco1--0} and Theorem \ref{t:wc}, we obtain the desired result. \qed

More concretely, we have the following example.

\begin{example}\label{exm1}
{\rm Let $\phi:\bR_+\to \bR_+$ be
a strictly increasing, continuous function such that  $\phi(0)=0$ and 
for all $0<r<R<\infty$, 
\[ c_1 \left(\frac Rr \right)^{\alpha_1} \leq\frac{\phi (R)}{\phi (r)} 
 \ \leq \ c_2 \left(\frac Rr \right)^{\alpha_2}
 \qquad \hbox{and} \qquad 
\int_0^r\frac {s}{\phi (s)}\, ds \leq  c_3 \, \frac{ r^2}{\phi
(r)}. \] 
Here $0<\alpha_1\le \alpha_2\le 2$. Assume that there exists
$\psi:\bR_+\to \bR_+$ a strictly increasing, continuous function with
$\psi(0)=0$ such that
\begin{equation}\label{eq:viwon}
\lim_{k\to\infty}\frac{\phi(k)}{\phi(kr)}=\frac 1{\psi(r)}
\qquad \hbox{for every } r>0.
\end{equation}

\noindent (i) Let $\{\xi_{xy}\}_{x,y\in \bZ^d, x\ne y}$ be
i.i.d. on $(\Omega, \sF, {\bf P})$ such that
$0\le\xi_{xy}$, $E[\xi_{xy}]=1$ and $\mbox{Var }(\xi_{xy})<\infty$.
Let $C(x,y)=\frac{\xi_{xy}}{|x-y|^{d}\phi(|x-y|)}$ for $x,y\in \bZ^d$, and define
instead of \eqref{eq:eneRCM},
\[
C^{(k)}(x,y):=k^{d} \phi(k) C(kx,ky)=\frac{\xi_{kx, ky} \phi(k)}{|x-y|^{d} \phi(k|x-y|) }
\qquad \mbox{for } x,y\in V_k.
\]
Then the claim of Proposition \ref{arrc}(i) holds, where
$X^{(k)}_t:=k^{-1}X^{(1)}_{\phi(k) t}$ and  $X$ is the Hunt process where the jump kernel
of the Dirichlet form is $j(x,y)=(|x-y|^{d}\psi(|x-y|))^{-1}$.

\noindent (ii) Assume further that $0\le\xi_{xy}\le C_1$ for some deterministic constant $C_1>0$. Then
the claim of Proposition \ref{arrc}(ii) holds.

\proof The proof of Proposition \ref{arrc} works line by line by plugging
$\frac{\phi(k)}{|x-y|^{d} \phi(k|x-y|) }$ into $j(x,y)$. Note that
instead of \eqref{eq:nirbw}, the following holds by using \eqref{eq:viwon},
\[  \lim_{k\to \infty}
k^{-2d}
\sum_{{x, y\in V_k\cap B_j}\atop{|x-y|>\delta}}
 {h(x,y)}\frac{\phi(k)}{|x-y|^{d} \phi(k|x-y|) }
= \int_{   B_j\times B_j} h(x,y) \, \frac{\1_{\{ |x-y|>\delta\}}}{|x-y|^{d} \psi(|x-y|) }m(dx)m(dy).\]
Given this equality, we can obtain {\bf (A4)${}^*$}\, (iii) by the same way as
that of Proposition \ref{arrc}. \qed }
\end{example}

\begin{remark} \label{r:6.2} \rm
(i) For the case of $d=1$, the only constraint is that the right hand side of
(\ref{eq:consrst}) is not summable. We can however obtain the corresponding results (strong
convergence of the semigroup and weak convergence) for
any subsequence $\{n_k\}$ such that $\sum_k1/n_k<\infty$.

\smallskip

(ii) The most typical case in the
Example \ref{exm1} is to take $\phi(r)=r^\alpha$.
Then $X^{(k)}_t=k^{-1}X^{(1)}_{k^\alpha t}$.
Thus Theorem \ref{arrc} says that,
if $d \ge 2$, $0\le\xi_{x,y}$,
${\bf E}[\xi_{x,y}] =1$ and $\mbox{Var }(\xi_{x,y})<\infty$, then
for any positive function $\varphi \in C_{c} (M)$,
  $\{(k^{-1}X^{(1)}_{k^\alpha t}, \, \P^{(k)}_\varphi); \, k \ge 1\}$
 converges weakly to $(X, \P_\varphi)$ on $\bD_{M_\partial} [0, 1]$ equipped with the convergence-in-measure
 topology ${\bf P}$-a.s., which in particular implies the finite dimensional convergence. Assume further
 that $0\le\xi_{x, y}\le C$
${\bf P}$-a.s.,
$\{(k^{-1}X^{(1)}_{k^\alpha t}, \, \P^{(k)}_\varphi); \, k
\ge 1\}$  converges weakly to $(X, \, \P_\varphi)$
on $\bD_{M_\partial} [0, 1]$ equipped with the Skorohod topology
${\bf P}$-a.s.. 

\smallskip

(iii) As mentioned in the introduction, one cannot obtain the a priori
H\"older estimates
 of caloric functions
in general (see \cite[Theorem 1.9]{BBCK}).  

\smallskip

(iv) It would be very nice if one can prove the Mosco convergence
for random walk on long range percolation.
Unfortunately, {\bf (A4)${}^*$}(ii) does not hold for the corresponding generator,
so we cannot apply Theorem \ref{Mosco1--0} for this model.
We note that the heat kernel bounds are obtained recently
in \cite{CS} for random walk on the long range percolation.
\end{remark}

\section{Appendix}

This appendix contains several equivalence conditions for
generalized Mosco convergence that was first obtained in
\cite[Theorem 2.5]{K} (appeared earlier in author's thesis
\cite{K0}).
In fact, a similar and more general form of such equivalence conditions for
generalized Mosco convergence was discussed in \cite{KS} independently.
Since we are using a minor modified version of \cite[Theorem 2.5]{K} and only the proof of (i) $\Longrightarrow$ (iv) is given in \cite{K}, we give full details for readers' convenience. We believe that, even if the version in \cite{KS} is quite general, our version in this paper is quite simple,
and it is
applicable to many cases.

For $k \ge 1$,   $(\sH_k, \langle\cdot, \cdot\rangle_k )$  and
$(\sH, \langle\cdot, \cdot\rangle)$ are Hilbert spaces
with the corresponding norms $\|\cdot\|_k$ and  $\|\cdot\|$.
  Suppose that $ (a^{(k)}, \sD(a^{(k)}))$
and $ (a, \sD(a))$ are densely defined closed symmetric
bilinear forms on $\sH^{(k)}$ and $\sH$, respectively.
We extend the definition of $a^{(k)}(u, u)$ to every $u\in \sH^{(k)}$
by defining $a^{(k)}(u,u) = \infty$ for   $u\in \sH_k \setminus \sD[a^{(k)}]$.
Similar extension is done for $a$ as well.

We assume throughout this section that for each $k\geq 1$,
 there is a  bounded linear operator $E_k : \sH_k \to \sH$ such that
 $\pi_k := E^*_k$ is a left inverse of $E_k$, that is,
\begin{equation}\label{eqn:AA}
 \langle \pi_k f, f_k \rangle_k = \langle f, E_k f_k \rangle
 \quad \hbox{and} \quad \pi_k E_k f_k = f_k  \qquad \hbox{for every } f \in  \sH, f_k \in \sH_k.
\end{equation}
Moreover we assume that $\pi_k: \sH \to \sH_k$ satisfies the following two conditions
\begin{equation}\label{eqn:AAA}
\sup_{k \ge 1} \|\pi_k  \|   \, < \infty,
\end{equation}
where $\| \pi_k\| $ denotes  the operator norm of $\pi_k$, and
\begin{equation}\label{eqn:A}
\lim_{k \to \infty} \|\pi_k f \|_k = \|f\| \qquad \hbox{for every } f \in \sH,
\end{equation}

Let $\| E_k\|$ denote the operator norm of $E_k: \sH^{(k)}\to \sH$.
Note that
$
\langle E_k f_k , E_k g_k \rangle = \langle f_k, g_k\rangle_k$  for every  $f_k, g_k \in \sH_k$, $k \ge 1$
and so clearly
\begin{equation}\label{eqn:A0}
\|E_k\| \equiv  1 ~~\mbox{ and }~~ \|E_k f_k\| = \|f_k\|_k ~~~\mbox{ for every } f_k \in \sH_k,   ~~k \ge 1
.
\end{equation}

\begin{definition}\label{D:Mosco}\rm
Under the above setting, we say that the closed bilinear form
$a^k$ is Mosco-convergent to $a$ in the generalized sense if
\begin{description}
\item{(i)}
If $v_k \in \sH_k$,  $u \in \sH$ and $E_k  v_k \to u$ weakly in $\sH$, then
$$
\liminf_{k \to \infty} a^{(k)} ( v_k, v_k ) \ge a(u,u) .
$$
\item{(ii)}
For every $u \in \sH$, there exists $u_k \in \sH_k$ such that
$f \in \sH$
$ E_k u_k  \to  u$ strongly in $\sH$
and
$$
\limsup_{k \to \infty} a^{(k)} (u_k, u_k) \le a(u,u).
$$
\end{description}
\end{definition}

Let $\big\{T^{(k)}_t, t\geq 1  \big\}$ and $\big\{G_{\lambda}^{(k)},
\lambda >0\big\}$   be the strongly continuous symmetric
contraction semigroup and the resolvent  associated with
$(a^{(k)}, \sD(a^{(k)}))$.  The infinitesimal generator of
$ \{T^{(k)}_t, t\geq 1\}$ (equivalently, of $(a^{(k)}, \sD(a^{(k)}))$)
will be denoted by
$A^k$.
Similarly, the semigroup, resolvent and
infinitesimal generator associated with $(a, \sD(a))$ will be denoted by
  $\{T_t, t\geq 0\}$,  $\left\{G_{\lambda}, \lambda>0\right\}$ and
  $A$  respectively.

\begin{theorem}\label{t:equi} Under the above setting,
 the followings are equivalent.
\begin{description}
\item{\rm (i)} $a^{(k)}$ is Mosco-convergent to $a$ in the generalized sense;

\item{\rm (ii)}
$E_k T^{(k)}_t \pi_k \to T_t$ strongly in $\sH$ and the convergence is uniform in any finite interval of $t \ge 0$;

\item{\em (iii)}
For each $f \in \sC_0$, there exists $\{f_k \}_{k \ge 1}$ such that $f_k \in \sD[\sA^{(k)}]$, $E_k f_k \to f$ and $E_k \sA^{(k)} f_k \to \sA f$ in $\sH$;

\item{(iv)}
$E_k G_{\lambda}^{(k)} \pi_k \to G_{\lambda}$ strongly in $\sH$ for every $\lambda > 0$.
\end{description}
\end{theorem}

\proof
Let $M_0 := \sup_{k  \ge 1} \|\pi_k\|$.
Note that,  by polarization identity and (\ref{eqn:A}), we have
\begin{equation}\label{eqn:A6}
\lim_{k \to \infty} \left\langle \pi_k u, \pi_k v \right\rangle_k = \langle u,v\rangle, \quad \text{for all } u, v \in \sH.
\end{equation}
By (\ref{eqn:AA})-(\ref{eqn:A0}), we see that for every $f \in \sH$ and $f_k \in \sH_k$,
\begin{eqnarray*}
\lim_{k \to \infty} \| f_k - \pi_k  f\|_k^2
&=&\lim_{k \to \infty}\left(\| f_k \|_k^2 -2 \langle f_k,  \pi_k f \rangle_k +\| \pi_k f\|_k^2\right)\\
&=&\lim_{k \to \infty} \left(\|E_k  f_k\|^2 -2\langle E_kf_k, f\rangle +  \| f\|^2 \right)
= \lim_{k \to \infty} \|E_k  f_k- f\|^2.
\end{eqnarray*}
Therefore
\begin{equation}\label{eqn:A2}
\lim_{k \to \infty} \left\|T^{(k)}_t \pi_k f - \pi_k T_t f\right\|_k
= \lim_{k \to \infty} \left\|E_k T^{(k)}_t \pi_k f - T_t f\right\|
\end{equation}
for every $f \in \sH$ and
\begin{equation}\label{eqn:A3}
\lim_{k \to \infty} \left\|G_{\lambda}^{(k)} \pi_k f - \pi_k G_{\lambda} f\right\|_k
= \lim_{k \to \infty} \left\|E_k G_{\lambda}^{(k)} \pi_k f - G_{\lambda} f \right\|
\end{equation}
for every $f \in \sH$ and $\lambda > 0$.

(ii) $\Longleftrightarrow$ (iii) : It is a special case of Theorem 1.6.1 in \cite{EK}.

(ii) $\Longleftrightarrow$ (iv) :
This can be proved using similar argument in the proof of Theorem 3.4.2 and Lemma 3.4.1 in \cite{P}.
We give a sketch here.
Similar to Lemma 3.4.1 in \cite{P}, one can check the following
\begin{equation}\label{eqn:A1_q}
E_k G_{\lambda}^{(k)} \left(\pi_k T_t - T^{(k)}_t \pi_k \right) G_{\lambda} f
~=~ \int^t_0 E_k T^{(k)}_{t-s} \left( \pi_k G_{\lambda}^{(k)} - G_{\lambda} \pi_k \right) T_s f ds
\end{equation}
for $f \in \sH$ and  $\lambda >0$. We first prove that (iv) implies (ii).

(ii) $\Longleftarrow$ (iv) :
We assume (iv) is true.
Fix  $\lambda >0$ and $T >0$,
If $f \in \sH$ and $ 0 \le t \le T$,
\begin{eqnarray*}
&& \left\| \left(E_k T^{(k)}_t \pi_k - T_t \right) G_{\lambda} f \right\| \\
&\le& \left\|E_k T^{(k)}_t \left(\pi_k G_{\lambda} - G_{\lambda}^{(k)} \pi_k \right)f \right\|
 + \left\| E_k G_{\lambda}^{(k)} \left(T^{(k)}_t \pi_k - \pi_k T (t)\right) f \right\|
 + \left\|\left( E_k G_{\lambda}^{(k)} \pi_k - G_{\lambda} \right) T_t f \right\| \\
&=& I_1 + I_2 + I_3 .
\end{eqnarray*}
$I_1 + I_3$ goes to $0$ uniformly on $[0,t]$ as $k \to \infty$ by (iv) and (\ref{eqn:A3}).
If $f \in \sD[A]$,
the domain of $A$,
there exists $g \in \sH$ such that $f = G_{\lambda} g$.
Since
$$
\left\|E_k T^{(k)}_{t-s}  \left(\pi_k G_{\lambda}T_s - G_{\lambda}^{(k)} \pi_k T_s \right)g\right\|
\le M_0 \left\|G_{\lambda} T_s g\right\| + \left\|G_{\lambda}^{(k)} \pi_k T_s g\right\|
\le \frac{2M_0}{\lambda} \|g\| ,
$$
by (\ref{eqn:A1_q}) and Lebesgue's dominated convergence theorem, we have
\begin{eqnarray*}
I_2 &\le& \int^t_0 \left\|E_k T^{(k)}_{t-s} \left(\pi_k G_{\lambda}T_s - G_{\lambda}^{(k)} \pi_k T_s \right) g\right\|ds \\
&\le& \int^t_0 \left\|\pi_k G_{\lambda} T_s g - G_{\lambda}^{(k)} \pi_k T_s g\right\|_k ds \rightarrow 0
\end{eqnarray*}
 uniformly on $[0,T]$ as $k \to \infty$ by (iv)
 and (\ref{eqn:A3}).
Since $A$ is densely defined, the above implies that (ii) is true.

(ii) $\Longrightarrow$ (iv):
Assume now that (ii) holds.
Then for $\lambda >0$ and $f \in \sH$,
$$
\left\|E_k G_{\lambda}^{(k)} \pi_k f - G_{\lambda} f \right\|
\,\le\, \int^{\infty}_0 e^{-\lambda t} \left\| \left(E_k T^{(k)}_t \pi_k - T_t \right)f \right\|dt \rightarrow 0~~~~ \mbox{ as } k \to \infty .
$$

(iv) $\Longrightarrow$ (i) :
Let
$$
a_{\lambda} (u, v) := \lambda\left\langle u - \lambda G_{\lambda} u,\, v \right\rangle ~~~~\mbox{ for } u, v \in \sH
$$
and
$$
a_{\lambda}^{(k)} (u_k, v_k) := \lambda\big\langle u_k - \lambda G_{\lambda}^k u_k, \, v_k\big\rangle_k ~~\mbox{ for } u_k, v_k \in \sH_k .
$$
It is well known that $a_{\lambda} (u, u)$ and $a_{\lambda}^{(k)} (u_k, u_k)$ are non-decreasing, and
$ \lim_{\lambda \to \infty} a_{\lambda} (u, u) =a (u, u)$ and
$ \lim_{\lambda \to \infty} a_{\lambda}^{(k)} (u_k, u_k)= a^k (u_k, u_k)$ for every $u \in \sH$ and $u_k\in \sH_k$.

Assume (iv) is true. By (\ref{eqn:A3}) and \eqref{eqn:A},
\begin{equation}\label{W_eqn:A4}
\lim_{k \to \infty} \left\| \left( G_{\lambda}^k \pi_k - \pi_k G_{\lambda} \right) f \right\|_k
= \lim_{k \to \infty} \left\| E_k G_{\lambda}^k \pi_k f - G_{\lambda} f \right\| = 0, \quad
\lim_{k \to \infty} \left\| G_{\lambda}^k \pi_k f \right\|_k = \|G_{\lambda}f \|
\end{equation}
for every $f \in \sH$ and $\lambda > 0$.
Since
\begin{eqnarray*}
&&|\,\lambda \langle \pi_k u - \lambda G_{\lambda}^k \pi_k u, \,\pi_k u \rangle_k - \lambda \langle u - \lambda G_{\lambda}u, \,u  \rangle \,| \\
&\le& \lambda^2 \,\|\, (G_{\lambda}^k \pi_k - \pi_k G_{\lambda} ) u \,\|_k\, \|\,\pi_k u\,\|_k
+\, \lambda\, |\, \langle \pi_k (u - \lambda G_{\lambda}u ),\, \pi_k u
\rangle_k - \langle u - \lambda G_{\lambda}u, \, u\rangle\, |,
\end{eqnarray*}
by (\ref{eqn:AAA}), (\ref{eqn:A6}) and (\ref{W_eqn:A4}) we have
\begin{equation}\label{W_eqn:A7}
\lim_{k \to \infty} a_{\lambda}^{(k)} (\pi_k u,\, \pi_k u)\, =\, a_{\lambda}(u,u) ~~~~\mbox{ for } \lambda >0.
\end{equation}
Suppose $v_k\in \sH_k$, $u \in \sH$ and $E_k v_k$ converges weakly to $u$ in $\sH$. By \eqref{eqn:AA} and \eqref{eqn:A6}
\begin{equation}\label{W_FFFF}
\lim_{k \to \infty} \left|\left\langle v_k - \pi_k u, \,\pi_k g\right\rangle_k\right| ~~~~\mbox{ for every } g \in \sH .
\end{equation}
We also have
$$
\lim_{k \to \infty} \langle v_k,\, \pi_k u\rangle_k\, =\, \|u\|^2,~~~~~ \sup_{k \ge 1} \| v_k\|_k < \infty
~~~~\mbox{ and }~~~~
\liminf_{k \to \infty} \| v_k\|_k ~\ge~ \|u\| .
$$
Note that
$$
a^k( v_k,  v_k) \ge a_{\lambda}^{(k)} ( v_k,~  v_k) \ge
 a_{\lambda}^{(k)} (\pi_k u, \pi_k u)  +2 \lambda \langle \pi_k u - \lambda G_{\lambda}^k \pi_k u,\,  v_k - \pi_k u \rangle_k.
$$
Since, by (iv) and
 (\ref{W_FFFF}),
\begin{eqnarray*}
 |\,\langle \pi_k u - \lambda G_{\lambda}^k \pi_k u,\,  v_k - \pi_k u \rangle_k \,|
&\le& |\, \langle \pi_k u, \, v_k - \pi_k u \rangle_k \,| \\
&&+ \lambda \,|\, \langle \pi_k G_{\lambda} u,\,  v_k - \pi_k u \rangle_k \,| \\
&&+ \lambda \, \|\, G_{\lambda}^k \pi_k u - \pi_k G_{\lambda}u \,\|_k \,( \|\, v_k\,\|_k \,+\, \|\,\pi_k u\,\|_k ),
\end{eqnarray*}
goes to $0$
 as $k \to \infty$ ,   we have by (\ref{W_eqn:A7}),
$$
\liminf_{k \to \infty}a^k( v_k,  v_k) \ge\liminf_{k \to \infty} a_{\lambda}^{(k)} ( v_k,~  v_k) \ge  \liminf_{k \to \infty}a_{\lambda}^{(k)} (\pi_k u, \pi_k u) = a_{\lambda}(u,u).
$$
Letting $\lambda \to \infty$, we obtain
$$
\liminf_{k \to \infty} a^k ( v_k, ~ v_k) ~\ge~ a(u,~u) .
$$

Now we suppose $u \in \sD[a]$ and show (ii) in Definition \ref{D:Mosco}.
First note that, by (iv),
$$
\lim_{\lambda \to \infty} \lambda \lim_{k \to \infty} E_k  G_{\lambda}^k \pi_k u =  \lim_{\lambda \to \infty}
\lambda   G_{\lambda} u = u, \quad \text{ in } \sH.$$
Thus, by (\ref{W_eqn:A7}) and the monotonicity of $a_{\lambda}^{(k)}$, we can choose an increasing sequence $\{ \lambda_k \}_{k \ge 1}$ such that
$$
\lim_{k \to \infty}\lambda_k = \infty, \quad \lim_{k \to \infty}   \lambda_k  E_k G_{\lambda_k}^k \pi_k u=u  \text{ in } \sH~~~~\mbox{ and }~~~~
\lim_{k \to \infty} a_{(\lambda_k)}^k ( \pi_k u,~ \pi_k u) \, \le\, a(u,\, u) < \infty .
$$
For $ k \ge 1 $, let
$
u_k := \lambda_k G_{\lambda_k}^k \pi_k u \in \sH_k$ and note that  $E_k u_k \to u$ in $\sH$.
Since
$$
a_{(\lambda_k)}^k ( \pi_k u, ~\pi_k u) = a^k (u_k, u_k) + \lambda_k \|u_k - \pi_k u\|^2_k =a^k (u_k, ~u_k) + \lambda_k \|E_k u_k - u\|^2,
$$
we conclude that
$$
a(u,\,u) \,\ge\, \limsup_{k \to \infty} a^k (u_k,~ u_k).
$$

(i) $\Longrightarrow$ (iv) :
Suppose (i) is true. Fix $\lambda > 0$ and assume $f \in \sH$.
Since
$$
\sup_{k \ge 1} \|E_k G_{\lambda}^{(k)} \pi_k \| \le \frac{M_0}{\lambda} < \infty,
$$
there exists a subsequence of $\left\{ E_k G_{\lambda}^{(k)} \pi_k f \right\}_{k \ge 1}$,
 still denoted $\left\{ E_k G_{\lambda}^{(k)} \pi_k f \right\}_{k \ge 1}$,
 such that $E_k G_{\lambda}^{(k)} \pi_k f $ converges weakly in $\sH$ to some $\tilde{u}$ in $\sH$.
So by Definition \ref{D:Mosco}(i)
\begin{equation}\label{AAA}
\liminf_{k \to \infty} \left(a^{(k)} (G_{\lambda}^{(k)} \pi_k f,\, G_{\lambda}^{(k)} \pi_k f )+
 \lambda\left\| G_{\lambda}^{(k)} \pi_k f \right\|^2_k \right)~\ge~
a( \tilde{u}, \tilde{u})+\lambda\left\| \tilde{u} \right\|^2.
\end{equation}
By \eqref{eqn:AA} and \eqref{AAA},
\begin{eqnarray}
&&a(  \tilde{u} ,  \tilde{u}) + \lambda \left\|   \tilde{u} \right\|^2 -2 \langle f,  \tilde{u}  \rangle  \nonumber \\
&\le& \liminf_{k \to \infty} \left( a^{(k)} ( G_{\lambda}^{(k)} \pi_k f , G_{\lambda}^{(k)} \pi_k f)
+ \lambda \left\| G_{\lambda}^{(k)} \pi_k f\right\|_k^2 \right) -2 \lim_{k \to \infty} \langle  f,E_k G_{\lambda}^{(k)} \pi_k f \rangle  \nonumber \\
&\le& \liminf_{k \to \infty} \left( a^{(k)} ( G_{\lambda}^{(k)} \pi_k f , G_{\lambda}^{(k)} \pi_k f)
+ \lambda \left\| G_{\lambda}^{(k)} \pi_k f\right\|_k^2 -2 \langle \pi_k f, G_{\lambda}^{(k)} \pi_k f \rangle_k \right) \\
&\le& \limsup_{k \to \infty} \left( a^{(k)} ( G_{\lambda}^{(k)} \pi_k f , G_{\lambda}^{(k)} \pi_k f)
+ \lambda \left\| G_{\lambda}^{(k)} \pi_k f\right\|_k^2 -2 \langle \pi_k f, G_{\lambda}^{(k)} \pi_k f \rangle_k \right)\label{e:AP}.
\end{eqnarray}

For arbitrary $v \in \sH$,  by Definition \ref{D:Mosco}(ii),
there exist $v_k \in \sH_k$ such that
\begin{equation}\label{e:B1}
\lim_{k \to \infty}\| v_k\|_k \,=\, \|v\|,~
 \lim_{k \to \infty} \langle E_k v_k, f \rangle =  \langle v, f \rangle
~\text{ and }~
\limsup_{k \to \infty} a^{(k)} (u_k, u_k) \le a(u,u).
\end{equation}

Since $G_{\lambda}^{(k)} \pi_k f$ is the unique minimizer of
$a^{(k)}(\,\cdot\, , \,\cdot\,) + \lambda \left\| \,\cdot\, \right\|_k^2 -2 \langle \pi_k f,\,\cdot\,  \rangle_k$
over $\sH_k$ for each $ k \ge 1$, \eqref{e:AP} is less than or equals to
$$
\limsup_{k \to \infty} a^{(k)} ( v_k , v_k)
+ \lambda \limsup_{k \to \infty} \left\| v_k \right\|_k^2 -2 \liminf_{k \to \infty}\langle \pi_k f, v_k \rangle_k,
$$
By (\ref{e:B1}), the above is less than or equals to
 $ a ( v ,  v) + \lambda \left\|   v \right\|^2 -2 \langle  f, v\rangle.$
Therefore $\tilde{u} = G_{\lambda} f$ because $G_{\lambda} f$ is  the unique minimizer of
$a(\,\cdot\, , \,\cdot\,) + \lambda \left\| \,\cdot\, \right\|^2 -2 \langle f,\,\cdot\,  \rangle$
over $\sH$.

On the other hand, by (i) there exists $w_k \in \sH_k$ such that
$$
\lim_{k \to \infty} \left\|E_kw_k -  G_{\lambda} f \right\| = 0
~~
\mbox{ and }
~~
\lim_{k \to \infty} a^{(k)}(w_k, ~w_k) = a(G_{\lambda} f,\, G_{\lambda} f) .
$$
So by (\ref{AAA}), the second equation above
 and the unique minimizer argument used above, we have
\begin{eqnarray*}
&&\lambda \limsup_{k \to \infty} \left\| G_{\lambda}^{(k)} \pi_k f - \frac{\pi_k f}{\lambda}\right\|^2_k \\
&\le& \limsup_{k \to \infty} \left(a^{(k)} (w_k, w_k) -  a^{(k)}(G_{\lambda}^{(k)} \pi_k f,~ G_{\lambda}^{(k)} \pi_k f)
+ \lambda \left\| w_k - \frac{\pi_k f}{\lambda}\right\|^2_k\right) \\
&\le& \limsup_{k \to \infty} a^{(k)} (w_k, w_k) - \liminf_{k \to \infty} a^{(k)}(G_{\lambda}^{(k)} \pi_k f,~ G_{\lambda}^{(k)} \pi_k f)
+ \lambda \limsup_{k \to \infty} \left\| w_k - \frac{\pi_k f}{\lambda}\right\|^2_k \\
&\le&
\lambda \limsup_{k \to \infty}
\left\| w_k - \frac{\pi_k f}{\lambda} \right\|^2_k .
\end{eqnarray*}
Combining the above inequality with
$$
\lim_{k \to \infty} \langle G_{\lambda}^{(k)} \pi_k f,\, \pi_k f \rangle_k
~=~ \langle G_{\lambda} f, f\rangle
~~
\mbox{ and }
~~
\lim_{k \to \infty} \left|\langle \pi_k f, w_k \rangle_k - \langle f, G_{\lambda} f
\rangle \right| ~=~ 0,
$$
we obtain
\begin{eqnarray*}
&&\limsup_{k \to \infty} \| G_{\lambda}^{(k)} \pi_k f \|_k\\
&=&\limsup_{k \to \infty}
 \left\| G_{\lambda}^{(k)} \pi_k f - \frac{\pi_k f}{\lambda}\right\|^2_k
+2 \lim_{k \to \infty} \langle G_{\lambda}^{(k)} \pi_k f,\, \pi_k f \rangle_k - \lim_{k \to \infty} \left\| \frac{\pi_k f}{\lambda} \right\|^2_k\\
&\le&\limsup_{k \to \infty}
 \left\| w_k - \frac{\pi_k f}{\lambda} \right\|^2_k
+2 \lim_{k \to \infty} \langle G_{\lambda}^{(k)} \pi_k f,\, \pi_k f \rangle_k - \lim_{k \to \infty} \left\| \frac{\pi_k f}{\lambda} \right\|^2_k\\
&\le&\limsup_{k \to \infty}
 \left\| w_k - \frac{\pi_k f}{\lambda} \right\|^2_k
+2 \lim_{k \to \infty} \langle w_k,\, \pi_k f \rangle_k - \lim_{k \to \infty} \left\| \frac{\pi_k f}{\lambda} \right\|^2_k\\
&=&\limsup_{k \to \infty}
\left\| (w_k - \frac{\pi_k f}{\lambda}) + \frac{\pi_k f}{\lambda} \right\|^2_k
\,=\, \limsup_{k \to \infty} \| w_k \|_k.
\end{eqnarray*}
Therefore
$$
\limsup_{k \to \infty} \| E_k G_{\lambda}^{(k)} \pi_k f \|
=\limsup_{k \to \infty} \| G_{\lambda}^{(k)} \pi_k f \|_k
~\le~ \limsup_{k \to \infty} \|E_k w_k \| = \|G_{\lambda} f \|
$$
and we conclude that,  for every $f \in \sH$, $E_k G_{\lambda}^{(k)} \pi_k f$ converges to $G_{\lambda} f$ in $\sH$.
\qed

\vskip 0.3truein

{\bf Zhen-Qing Chen}

Department of Mathematics, University of Washington, Seattle,
WA 98195, USA

E-mail: zchen@math.washington.edu

\bigskip

{\bf Panki Kim}

Department of Mathematical Sciences and Research Institute of Mathematics,

Seoul National University,
San56-1 Shinrim-dong Kwanak-gu,
Seoul 151-747, Republic of Korea

E-mail: pkim@snu.ac.kr

\bigskip

{\bf Takashi Kumagai}

Department of Mathematics, Faculty of Science,
Kyoto University, Kyoto 606-8502, Japan

E-mail: kumagai@math.kyoto-u.ac.jp

\vfill
\end{document}